%% file: article2.tex
\documentclass[12pt,dvips]{amsart}
\usepackage[usenames]{color}
\usepackage{tikz}

\newtheorem{lemma}{Lemma}[section]
\newtheorem{theorem}{Theorem}[section]
\newtheorem{definition}{Definition}[section]
\newtheorem{proposition}{Proposition}[section]

\DeclareMathOperator{\bin}{Bin}

\DeclareMathOperator{\CPP}{CPP}
\DeclareMathOperator{\RPP}{RPP}

\DeclareMathOperator{\lhs}{lhs}

\DeclareMathOperator{\PP}{PP}

\DeclareMathOperator{\lab}{label}

\DeclareMathOperator{\ALCD}{ALCD}
\DeclareMathOperator{\IP}{P}

\DeclareMathOperator{\hook}{hook}
\DeclareMathOperator{\cohook}{cohook}

\DeclareMathOperator{\diag}{diag}

\DeclareMathOperator{\cols}{columns}

\title{Enumeration of Cylindric Plane Partitions - part II}
\author{Robin Langer}
\date{\today}

\begin{document}
\input tableau.tex
%
\mbox{}
\vspace{2ex}

\maketitle

\begin{abstract}
Cylindric plane partitions may be thought of as a natural generalization of reverse plane partitions. A generating series for the enumeration of cylindric plane partitions was recently given by Borodin. As in the reverse plane partition case, the right hand side of this identity admits a simple factorization form in terms of the ``hook lengths'' of the individual boxes in the underlying shape. The main result of this paper is a new bijective proof of Borodin's identity which makes use of Fomin's growth diagram framework for generalized RSK correspondences.
\end{abstract}

\section{Introduction}
\label{sec:in}
Cylindric plane partitions were first introduced by Gessel and Krattenthaler \cite{gessel-1997}. 
For any binary string $\pi$ of length $T$,
a \emph{cylindric plane partition} with profile $\pi$ may be defined as a sequence of integer partitions:
\begin{equation} 
(\mu^0, \mu^1, \ldots \mu^T) \qquad \qquad \mu^0 = \mu^T 
\end{equation}
such that if $\pi_k = 1$ then $\mu^k / \mu^{k-1}$ is a \emph{horizontal strip}.
Otherwise if $\pi_k = 0$ then $\mu^{k-1} / \mu^k$ is a horizontal strip.
The \emph{weight} of a cylindric partition is given by
$  |\mathfrak{c}| = |\mu_1| + |\mu_2| + \cdots |\mu_T|$.
In the special case where $\mu_0 = \mu_T = \emptyset$ we recover the usual definition of a reverse plane partition \cite{adachi}.
If, in addition to this there are no inversions in the profile, we have a regular plane partition.

For those readers who are more familliar with the definition of a plane partition as an array of integers which is weakly decreasing along both rows and columns, the bijection with the ``interlacing sequence'' model is obtained by reading along the main diagonals. For example:
\[  \tableau{3 & 3 & 2 \\ 3 & 2 & 1 \\ 1 & 1 & 1} \]
\[ \mathfrak{c} = (\emptyset, (1), (3,1), (3,2,1), (3,1), (2), \emptyset ) \]

A regular plane partition may also be thought of as a pair of semi-standard young tableau of the same shape. In the case of our example, the two tableaux are:
\[
 \tableau{1 & 2 & 2 \\ 2 & 3 \\ 3} \qquad \tableau{1 & 1 & 2 \\ 2 & 3  \\ 3}
\] 

The theory of plane partitions is closely related to both the theory of symmetric functions and Fomin's theory of generalized RSK type correspondences \cite{fomin1,fomin2}.
The beginning of the subject is perhaps the following famous identity of MacMahon:
\begin{equation} \label{macmahon}
\sum_{\mathfrak{c} \in \PP} z^{|\mathfrak{c}|} = \left ( \frac{1}{1-z^n} \right )^n 
\end{equation}

It was first pointed out by Okounkov \cite{okounkov} that enumerative results for plane partitions may be obtained by considering commutation relations between vertex operators acting on fermionic fock space. The underlying algebraic structure is that of the Heisenberg algebra. By the boson-fermion correspondence these operators may be alternatively thought of as acting on symmetric functions.
The \emph{Pieri rules} for Schur functions are key to this approach:
\begin{equation}\label{p1} S_\mu[X] h_r[X] = \sum_{\lambda \in U_r(\mu)} S_\lambda[X] \end{equation}
\begin{equation}\label{p2} S_\lambda[X+z] = \sum_{\mu \in D_r(\mu)} S_\mu[X] z^r \end{equation}
Here $U_r(\mu)$ denotes the set of all partitions which can be obtained from $\mu$ by adding a horrizontal $r$-strip and $D_r(\lambda$) denotes the set of all partitions which can be obtained from $\lambda$ by removing a horizontal $r$-strip.

The next important result in the subject is the following hook-product formula for the enumeration of reverse plane partitions with arbitrary profile $\pi$ which is due to Stanley:
\begin{equation} \label{stanley}
\sum_{\mathfrak{c} \in \RPP(\pi)} z^{|\mathfrak{c}|} = 
\prod_{\substack{i < j \\ \pi_i > \pi_j}} \frac{1}{1-z^{j-i}} 
\end{equation}

There is a natural bijection between the boxes of a partition and the inversions in its profile. The expression $j-i$ on the right hand side of Stanley's identity may be understood as the \emph{hook length} of the corresponding box (see section \ref{sec:def}). 

In fact, the right hand side of Stanley's identity may be interpreted as a weighted sum over arbitrarily labelled young diagrams. The weight of such a labelled diagram is equal to a sum over the boxes of the diagram, of the label of the box, times the hook length of the box. 
The Hillman-Grassl algorithm provides a bijective proof of Stanley's identity. 

More recently the following hook-product formula for the enumeration of cylindric plane partitions of given profile was first given by Borodin \cite{borodin}. A very different proof involving the representation theory of $\widehat{sl}(n)$ was later given by Tingley \cite{tingley}: 
\begin{equation} \label{borodin}
\sum_{\mathfrak{c} \in \CPP(\pi)} z^{|\mathfrak{c}|} = 
\prod_{n \geq 0} \left ( \frac{1}{1-z^{nT}} \prod_{\substack{i < j \\ \pi_i > \pi_j}} \frac{1}{1-z^{j-i + nT}} 
\prod_{\substack{i > j \\ \pi_i > \pi_j }} \frac{1}{1 - z^{j-i + (n+1)T}} \right ) 
\end{equation}

Here $T$ denotes the length of the profile. As in the reverse plane partition case, there is a natural bijection between the ``boxes'' of the cylindric plane partition and ``cylindric inversions'' of the ``underlying shape''. The expression $j-i + kT$ on the right hand side of Borodin's identity may be understood as the ``cylindric hook length'' of the box with ``cylindric inversion coordinates'' $(i,j,k)$ (see section \ref{sec:def}). 

The right hand side of Borodin's identity  may be interpreted combinatorially as a sum over pairs $(\gamma, \mathfrak{A})$ where $\gamma$ is an integer partition and $\mathfrak{A}$ is an arbitrarily labelled ``cylindric diagram''. The weight of $\mathfrak{A}$ is a sum over the boxes of $\mathfrak{A}$ of the label of the box times the cylindric hook length of the box. The weight of the pair $(\gamma, \mathfrak{A})$ is $T|\gamma| + |\mathfrak{A}|$.

The main result of this paper is a bijective proof of Borodin's identity:
\begin{theorem}\label{borodin-bijective}.
There exists a natural weight preserving bijection between cylindric plane partitions with given profile $\pi$ and pairs $(\gamma, \mathfrak{A})$ where $\gamma$ is an integer partition and $\mathfrak{A}$ is an arbitrarily labelled cylindric diagram with profile $\pi$.
\end{theorem}
The proof uses the idea of local rules and growth diagrams first introduced by Fomin \cite{fomin1, fomin2}. The reverse plane partition version of this proof has been previously given by Krattenthaler \cite{krattenthaler}.

\section{Definitions}
\label{sec:def}

\subsection{Partitions and profiles}\label{profiles}

An \emph{integer partition} is simply a weakly decreasing list of non-negative integers which eventually stabilizes at zero. 
If the sum of the parts of $\lambda$ is equal to $n$, then we say that $\lambda$ is a partition of $n$ and write $|\lambda| = n$. 
The conjugate of the integer partition $\lambda = (\lambda_1, \lambda_2, \ldots, \lambda_k)$ is defined to be $\lambda' = (\lambda'_1, \lambda'_2, \ldots \lambda'_r)$ where $\lambda'_j = \#\{i \,|\, \lambda_i \geq j \}$.

It is often convenient to represent an integer partition visually as a \emph{Young diagram}, which is a collection of boxes in the cartesian plane which are ``stacked up'' in the bottom right hand corner. Note that our convention differs from both the standard French and English conventions. 

The \emph{profile} of an integer partition is the binary string which traces out the ``jagged boundary'' of the associated young diagram. Reading from the top right hand corner to the bottom left hand corner, a zero is recorded for every vertical step and a one for every horizontal step. For example the minimum profile of our example partition $\lambda = (5,3,3,2)$ is $110100110$:

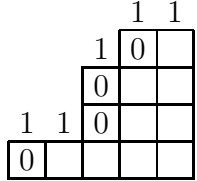
\[
\tableau{
\missingcell & \missingcell & \missingcell & \missingcell 1 & \missingcell 1 \\
\missingcell & \missingcell & \missingcell 1 & 0 & \\
\missingcell & \missingcell & 0 & &  \\
\missingcell 1 & \missingcell1 & 0 & &  \\
0 & &&&\\
} 
\]

\begin{definition}\label{inversion}
An \emph{inversion} in a binary string $\pi$ is a pair of indices $(i,j)$ such that $i < j$ and $\pi_i > \pi_j$. 
\end{definition}
There is a natural bijection between the ``boxes'' of an integer partition $\lambda$ and the inversions in any generalized profile of $\lambda$. If the box $s$ has ``inversion coordinates'' $(i,j)$ then the hook length is given by $h_\lambda(s) = j - i$.

We shall make extensive use of the partial order on the set of all binary strings $\bin(n,m)$ with $n$ zeros and $m$ ones, whose covering relations are given by
$\pi \prec \pi'$ if and only if there is some $i$ such that $\pi_i = 0 = \pi'_{i+1}$ and $\pi_{i+1} = 1 = \pi'_i$
and for all other $k$ we have $\pi'_k = \pi_k$. 
In other words $\pi'$ is obtained from $\pi$ by adding an inversion.

We shall denote by $\pi_{\min}$ the binary string with $n$ zeros follows by $m$ ones
and $\pi_{\max}$ the binary string with $m$ ones follows by $n$ zeros.

\subsection{Cylindric Diagrams}\label{diagrams}

A \emph{cylindric diagram} may be thought of as an infinite partition with periodic profile, which has been wrapped around a cylinder. We shall use the notation $\hat{\lambda}(\pi)$ to denote the cylindric diagram with profile $\pi$.

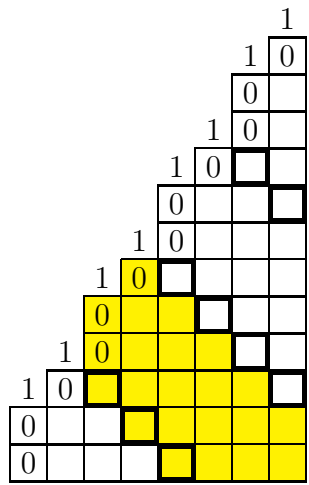
\[
 \tableau{
\missingcell & \missingcell & \missingcell & \missingcell & \missingcell & \missingcell & \missingcell & \missingcell 1  \\
\missingcell & \missingcell & \missingcell & \missingcell & \missingcell & \missingcell & \missingcell 1&0  \\
\missingcell & \missingcell & \missingcell & \missingcell & \missingcell & \missingcell &0& \\
\missingcell & \missingcell & \missingcell & \missingcell & \missingcell & \missingcell 1&0& \\
\missingcell & \missingcell & \missingcell & \missingcell & \missingcell 1&0& \thickcell&  \\
\missingcell & \missingcell & \missingcell & \missingcell &0&&  & \thickcell\\
\missingcell & \missingcell & \missingcell & \missingcell 1&0&&  & \\
\missingcell & \missingcell & \missingcell 1&\yellowcell 0& \thickcell &  && \\
\missingcell & \missingcell &\yellowcell0&\yellowcell&\yellowcell & \thickcell&& \\
\missingcell & \missingcell 1&\yellowcell0&\yellowcell& \yellowcell &\yellowcell&\thickcell& \\
\missingcell1 &0 & \thickyellowcell & \yellowcell &\yellowcell&\yellowcell&\yellowcell& \thickcell\\
0&&& \thickyellowcell & \yellowcell &\yellowcell&\yellowcell& \yellowcell\\
0&&&& \thickyellowcell &\yellowcell&\yellowcell& \yellowcell\\
}
\]
Profile $\pi = 10100$. Period $T = 5$. The ``fundamental domain'' is coloured in yellow. Only a finite section of the cylindric diagram is shown.

It is natural to index the boxes of the cylindric diagram via ``cylindric inversion coordinates'' $(i,j,k)$ where $\pi_i = 1$, $\pi_j = 0$ and if $j < i$ then $k \geq 1$ otherwise $k \geq 0$.
Here are the cylindric inversion coordinates of each box of our example cylindric diagram:

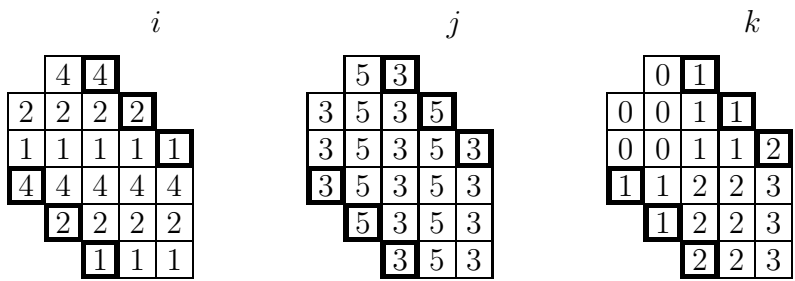
\begin{align*}
i \quad & \quad & j \quad & \quad & k \quad \\
\tableau{
\missingcell & 4 & \thickcell{4} \\
2 & 2 & 2 & \thickcell{2} \\
1 & 1 & 1 & 1 & \thickcell{1} \\
\thickcell{4} & 4 & 4 & 4 & 4 \\
\missingcell & \thickcell{2} & 2 & 2 & 2 \\
\missingcell & \missingcell & \thickcell{1} & 1 & 1 
}
& \quad &
\tableau{
\missingcell & 5 & \thickcell{3} \\
3 & 5 & 3 & \thickcell{5} \\
3 & 5 & 3 & 5 & \thickcell{3} \\
\thickcell{3} & 5 & 3 & 5 & 3 \\
\missingcell & \thickcell{5} & 3 & 5 & 3 \\
\missingcell & \missingcell & \thickcell{3} & 5 & 3 
}
& \quad &
\tableau{
\missingcell & 0 & \thickcell{1} \\
0 & 0 & 1 & \thickcell{1} \\
0 & 0 & 1 & 1 & \thickcell{2} \\
\thickcell{1} & 1 & 2 & 2 & 3 \\
\missingcell & \thickcell{1} & 2 & 2 & 3 \\
\missingcell & \missingcell & \thickcell{2} & 2 & 3 
} 
\end{align*}

Two boxes lie in the same ``cylindric row'' if they have the same $i$-coordinate, and in the same ``cylindric column'' if they have the same $j$-coordinate. 

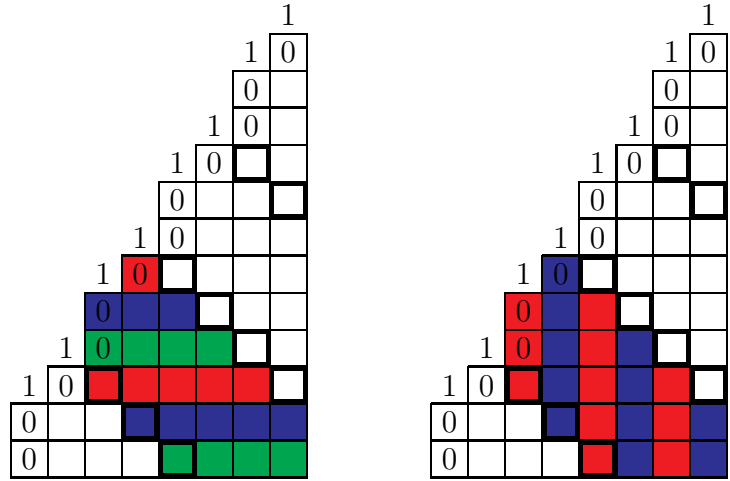
\[
 \tableau{
\missingcell & \missingcell & \missingcell & \missingcell & \missingcell & \missingcell & \missingcell & \missingcell 1  \\
\missingcell & \missingcell & \missingcell & \missingcell & \missingcell & \missingcell & \missingcell 1&0  \\
\missingcell & \missingcell & \missingcell & \missingcell & \missingcell & \missingcell &0& \\
\missingcell & \missingcell & \missingcell & \missingcell & \missingcell & \missingcell 1&0& \\
\missingcell & \missingcell & \missingcell & \missingcell & \missingcell 1&0& \thickcell&  \\
\missingcell & \missingcell & \missingcell & \missingcell &0&&  & \thickcell\\
\missingcell & \missingcell & \missingcell & \missingcell 1&0&&  & \\
\missingcell & \missingcell & \missingcell 1&\redcell0& \thickcell &  && \\
\missingcell & \missingcell &\bluecell0&\bluecell&\bluecell & \thickcell&& \\
\missingcell & \missingcell 1&\greencell0&\greencell& \greencell &\greencell&\thickcell& \\
\missingcell1 &0 & \thickredcell & \redcell &\redcell&\redcell&\redcell& \thickcell\\
0&&& \thickbluecell & \bluecell &\bluecell&\bluecell& \bluecell\\
0&&&& \thickgreencell &\greencell&\greencell&\greencell \\
}
\qquad \qquad
 \tableau{
\missingcell & \missingcell & \missingcell & \missingcell & \missingcell & \missingcell & \missingcell & \missingcell 1  \\
\missingcell & \missingcell & \missingcell & \missingcell & \missingcell & \missingcell & \missingcell 1&0  \\
\missingcell & \missingcell & \missingcell & \missingcell & \missingcell & \missingcell &0& \\
\missingcell & \missingcell & \missingcell & \missingcell & \missingcell & \missingcell 1&0& \\
\missingcell & \missingcell & \missingcell & \missingcell & \missingcell 1&0& \thickcell&  \\
\missingcell & \missingcell & \missingcell & \missingcell &0&&  & \thickcell\\
\missingcell & \missingcell & \missingcell & \missingcell 1&0&&  & \\
\missingcell & \missingcell & \missingcell 1& \bluecell 0& \thickcell &  && \\
\missingcell & \missingcell &\redcell0&\bluecell& \redcell & \thickcell&& \\
\missingcell & \missingcell 1&\redcell0&\bluecell&\redcell  &\bluecell&\thickcell& \\
\missingcell1 &0 & \thickredcell & \bluecell &\redcell&\bluecell&\redcell& \thickcell\\
0&&& \thickbluecell & \redcell & \bluecell & \redcell & \bluecell \\
0&&&& \thickredcell & \bluecell & \redcell & \bluecell \\
}
\]

The $k$-coordinate may be thought of as a sort of ``depth'' or ``winding number''.

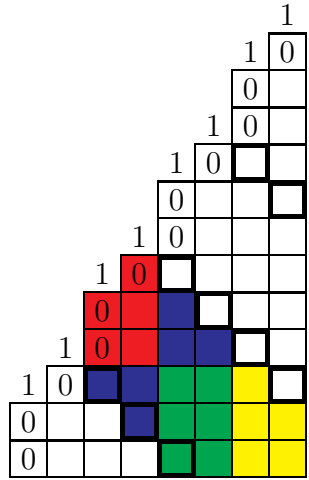
\[
 \tableau{
\missingcell & \missingcell & \missingcell & \missingcell & \missingcell & \missingcell & \missingcell & \missingcell 1  \\
\missingcell & \missingcell & \missingcell & \missingcell & \missingcell & \missingcell & \missingcell 1&0  \\
\missingcell & \missingcell & \missingcell & \missingcell & \missingcell & \missingcell &0& \\
\missingcell & \missingcell & \missingcell & \missingcell & \missingcell & \missingcell 1&0& \\
\missingcell & \missingcell & \missingcell & \missingcell & \missingcell 1&0& \thickcell&  \\
\missingcell & \missingcell & \missingcell & \missingcell &0&&  & \thickcell\\
\missingcell & \missingcell & \missingcell & \missingcell 1&0&&  & \\
\missingcell & \missingcell & \missingcell 1& \redcell 0& \thickcell &  && \\
\missingcell & \missingcell &\redcell0&\redcell& \bluecell & \thickcell&& \\
\missingcell & \missingcell 1&\redcell0&\redcell&\bluecell  &\bluecell&\thickcell& \\
\missingcell1 &0 & \thickbluecell & \bluecell &\greencell&\greencell&\yellowcell& \thickcell\\
0&&& \thickbluecell & \greencell & \greencell & \yellowcell & \yellowcell \\
0&&&& \thickgreencell & \greencell & \yellowcell & \yellowcell \\
} \\ \\
\]

The \emph{cylindric hook length} of a box is the hook length of the box relative to the larger partition.

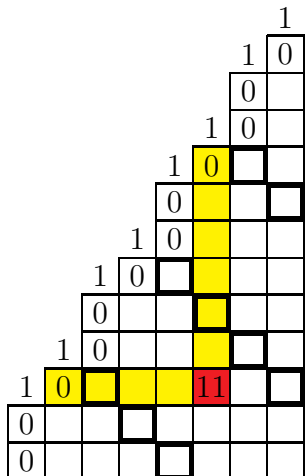
\[
 \tableau{
\missingcell & \missingcell & \missingcell & \missingcell & \missingcell & \missingcell & \missingcell & \missingcell 1  \\
\missingcell & \missingcell & \missingcell & \missingcell & \missingcell & \missingcell & \missingcell 1&0  \\
\missingcell & \missingcell & \missingcell & \missingcell & \missingcell & \missingcell &0& \\
\missingcell & \missingcell & \missingcell & \missingcell & \missingcell & \missingcell 1&0& \\
\missingcell & \missingcell & \missingcell & \missingcell & \missingcell 1&\yellowcell0& \thickcell&  \\
\missingcell & \missingcell & \missingcell & \missingcell &0&\yellowcell&  & \thickcell\\
\missingcell & \missingcell & \missingcell & \missingcell 1&0&\yellowcell&  & \\
\missingcell & \missingcell & \missingcell 1&0& \thickcell &\yellowcell  && \\
\missingcell & \missingcell &0&& & \thickyellowcell&& \\
\missingcell & \missingcell 1&0&&  &\yellowcell&\thickcell& \\
\missingcell1 & \yellowcell 0 &  \thickyellowcell & \yellowcell &\yellowcell&\redcell 11&& \thickcell\\
0&&& \thickcell &  &&& \\
0&&&& \thickcell &&& \\
}
\]

We shall use the notation $h_{\hat{\lambda}(\pi)}(b)$ to denote the cylindric hook length of the box $b$ relative to the cylindric diagram $\hat{\lambda}(\pi)$

\begin{lemma}\label{cylindrichook}
The cylindric hook length of a box with cylindric inversion coordinates $(i,j,k)$ is given by $j - i + kT$.
\end{lemma}

Here are the hook lengths of the boxes in our example partition:

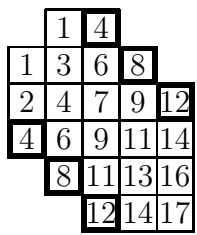
\[  \tableau{
\missingcell & 1 & \thickcell{4} \\
1 & 3 & 6 & \thickcell{8} \\
2 & 4 & 7 & 9 & \thickcell{12} \\
\thickcell{4} & 6 & 9 & 11 & 14 \\
\missingcell & \thickcell{8} & 11 & 13 & 16 \\
\missingcell & \missingcell & \thickcell{12} & 14 & 17 
}
\]

\begin{definition}\label{shift}
The rotation operator on binary strings is defined by:
\[ \sigma(\pi)_i = \pi_{(i + 1) \mod T}\]
where $T$ is the length of $\pi$.
\end{definition}

There is a natural bijection between the cylindric diagram with profile $\pi$ and the 
cylindric diagram with profile $\sigma(\pi)$. Nevertheless since the same box will have
different cylindric inversion coordinates, depending on the choice of rotation of the profile, 
we prefer to consider two cylindric diagrams which differ by a rotation to be two distinct objects.

\subsection{Cylindric Plane Partitions}\label{cpp}

For any binary string $\pi$ of length $T$,
a \emph{cylindric plane partition} with profile $\pi$ may be defined as a sequence of integer partitions:
\begin{equation} 
\mathfrak{c} = (\mu^0, \mu^1, \ldots \mu^T) \qquad \qquad \mu^0 = \mu^T 
\end{equation}
such that if $\pi_k = 1$ then $\mu^k / \mu^{k-1}$ is a \emph{horizontal strip}.
Otherwise if $\pi_k = 0$ then $\mu^{k-1} / \mu^k$ is a horizontal strip.

The \emph{weight} of a cylindric plane partition is given by:
\[  |\mathfrak{c}| = |\mu_1| + |\mu_2| + \cdots |\mu_T| \]
The \emph{refined weight} of a cylindric plane partition is given by:
\[  z_1^{|\mu_1|} z_2^{|\mu_2|}  \cdots z_T^{|\mu_T|} \]

In the special case where $\mu_0 = \mu_T = \emptyset$ we recover the usual definition of a reverse plane partition \cite{adachi}.
If, in addition to this there are no inversions in the profile, we have a regular plane partition.

\begin{lemma}
The rotation operator $\sigma$ naturally induces a weight preserving map:
\[ \sigma : \CPP(\pi) \to \CPP(\sigma(\pi)) \]
\end{lemma}

A \emph{cube} of a cylindric partition $\mathfrak{c} = (\mu_0, \mu_1, \ldots \mu_T)$ is a \emph{box} of any of the partitions $\mu_1, \mu_2, \ldots \mu_T$. 
Note that to avoid double counting we do not include the boxes in the partition $\mu_0$. 

Cylindric plane partitions are often represented as certain labelled cylindric diagrams. For example, the cylindric plane partition
\[ \mathfrak{c} = ((3,2,2), (5,3,2), (6,4,3,2), (4,3,2), (4,3,2,1), (3,2,2))\] 
with profile $10100$
may be represented as:
\[  \tableau{
\missingcell & 4 & \thickcell{3} \\
6 & 4 & 3 & \thickcell{2} \\
5 & 4 & 3 & 2 & \thickcell{2} \\
\thickcell{3} & 3 & 3 & 2 & 1 \\
\missingcell & \thickcell{2} & 2 & 2 & 0 \\
\missingcell & \missingcell & \thickcell{2} & 0 & 0 
}
\]
The individual partitions in the interlacing sequence picture are read off the diagonals.

\begin{lemma}
The labels of the cylindric diagram associated to a cylindric plane partition are weakly decreasing along both cylindric rows and cylindric columns.
\begin{proof}
This is an immediate consequence of the horizontal strip condition on diagonals
\end{proof}
\end{lemma}

\subsection{Arbitrarily Labelled Cylindric Diagrams}\label{alcd}

An \emph{arbitrarily labelled cylindric diagram} $\mathfrak{d}$ with profile $\pi$ is simply an assignment of non-negative integers to the boxes of the associated cylindrical diagram
in such a way that only finitely many of the labels are non-zero. We shall use the notation $\ALCD(\pi)$ to denote the set of all arbitrarily labelled cylindric diagrams with profile $\pi$.

\begin{definition} \label{depth}
The \emph{depth} of an arbitrarily labelled cylindric diagram $\mathfrak{d}$ is the smallest $k$ such that all boxes with cylindric inversion coordinates $(i,j,k')$ with $k' \geq k$ have label zero.
\end{definition}

The weight of an arbitrarily labelled cylindric diagram is given by the sum over boxes in the cylindric diagram of the label of the box times the hook length of the box:
\[ |\mathfrak{d}| = \sum_{b \in \hat{\lambda}(\pi)} \lab(b) \,\, |\hook(b)| \]

For example, the following arbitrarily labelled cylindric diagram has depth 2 and weight 26.

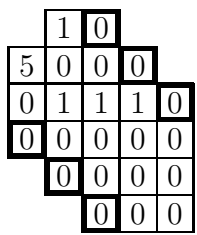
\[ \tableau{
\missingcell & 1 & \thickcell{0} \\
5 & 0 & 0 & \thickcell{0} \\
0 & 1 & 1 & 1 & \thickcell{0} \\
\thickcell{0} & 0 & 0 & 0 & 0 \\
\missingcell & \thickcell{0} & 0 & 0 & 0 \\
\missingcell & \missingcell & \thickcell{0} & 0 & 0 
}
\]

\begin{lemma}
The rotation operator $\sigma$ naturally induces a well-defined weight preserving map:
\[ \sigma : \ALCD(\pi) \to \ALCD(\sigma(\pi)) \]
\end{lemma}

\subsubsection{alternative definition of weight}

Let us define the coohook of a box $b$ in an arbitrarily labelled cylindric diagram $\mathfrak{d}$ to be:
\[ \cohook(b) = \{ b' | b \in \hook(b) \} \]

\[
 \tableau{
\missingcell & \missingcell & \missingcell & \missingcell & \missingcell & \missingcell & \missingcell & \missingcell 1  \\
\missingcell & \missingcell & \missingcell & \missingcell & \missingcell & \missingcell & \missingcell 1&0  \\
\missingcell & \missingcell & \missingcell & \missingcell & \missingcell & \missingcell &0& \\
\missingcell & \missingcell & \missingcell & \missingcell & \missingcell & \missingcell 1&0& \\
\missingcell & \missingcell & \missingcell & \missingcell & \missingcell 1&0& \thickcell&  \\
\missingcell & \missingcell & \missingcell & \missingcell &0&&  & \thickcell\\
\missingcell & \missingcell & \missingcell & \missingcell 1&0&&  & \\
\missingcell & \missingcell & \missingcell 1&0& \thickcell &  && \\
\missingcell & \missingcell &0& \redcell & \yellowcell & \thickyellowcell& \yellowcell& \yellowcell \\
\missingcell & \missingcell 1&0& \yellowcell &  &&\thickcell& \\
\missingcell1 &  0 &  \thickcell & \yellowcell && && \thickcell\\
0&&& \thickyellowcell &  &&& \\
0&&& \yellowcell & \thickcell &&& \\
}
\]

Although number of boxes in the cohook of a given box is always infinite, for any given arbitrarily labelled cylindric plane partition it is meaningful to define the weight of a cohook:
\[ |\cohook(b)|_\mathfrak{d} = \sum_{b' \in \cohook(b)} \lab(b') \]

Let $\diag(k)$ denote the set of all boxes on the $k$-th diagonal.
Furthermore let us define the weight of the diagonal of an arbitrarily labelled cylindric diagram $\mathfrak{d}$ to be:
\[ |\diag(k)|_\mathfrak{d} = \sum_{b \in \diag(k)} |\cohook(b)|_\mathfrak{d} \]

With these definitions we may give an alternative definition of the weight of $\mathfrak{d}$:
\begin{align*}
|\mathfrak{d} |
& = \sum_{b' \in \lambda} \lab(b') \hook(b)  \\
& = \sum_{b' \in \lambda} \lab(b') \sum_{b \in \hook(b')} 1  \\
& =  \sum_{b \in \lambda} \, \sum_{b' \in \cohook(b)} \lab(b') \\
& = \sum_{b \in \lambda} |\cohook(b)|_\mathfrak{d} \\
& = \sum_{k=1}^T |\diag(k)|_\mathfrak{d} \\
\end{align*}
The \emph{refined weight} of a cylindric plane partition is given by:
\[  z_1^{|\diag(1)|_\mathfrak{d}} z_2^{|\diag(2)|_\mathfrak{d}} \cdots z_T^{|\diag(T)|_\mathfrak{d}} \]

Note that if $b$ is taken to be the box of $\diag(k)$ lying furthest to the ``north-west'' then $|\diag(k)|_\mathfrak{d}$ is none other than
the sum of all the labels of boxes lying ``south-east'' of $b$.

\section{Combinatorial Interpreation of Borodin's Identity}
\label{statement}

Borodin's identity for the enumeration of cylindric plane partitions \cite{borodin} states that for a given profile $\pi$ of length $T$ we have:
\begin{equation} \label{borodin}
\sum_{\mathfrak{c} \in \CPP(\pi)} z^{|\mathfrak{c}|} = 
\prod_{n \geq 0} \left ( \frac{1}{1-z^T} \prod_{\substack{i < j \\ \pi_i > \pi_j}} \frac{1}{1-z^{j-i + nT}} 
\prod_{\substack{i < j \\ \pi_i < \pi_j }} \frac{1}{1 - z^{j-i + (n+1)T}} \right ) 
\end{equation}
Here $\CPP(\pi)$ denotes the set of all cylindric plane partitions with profile $\pi$ (see section \ref{cpp}).

As a consequence of Lemma \ref{cylindrichook} Borodin's identity may be rewritten in the form:

\begin{align*}
\sum_{\mathfrak{c} \in \CPP(\pi)} z^{|\mathfrak{c}|} & = \left ( \sum_{\gamma} \frac{1}{1-z^{T\, |\gamma|}} \right ) 
\left ( \sum_{s \in \widehat{\lambda}(\pi)} \frac{1}{1 - z^{h_{\widehat{\lambda}(\pi)}(s)}} \right ) 
\end{align*}
where $\widehat{\lambda}(\pi)$ denotes the \emph{cylindric diagram} with profile $\pi$
and $h_{\widehat{\lambda}(\pi)}(s)$ denotes the \emph{cylindric hook length} of the box $s$ (see Section \ref{diagrams}).

The right hand side may be interpreted combinatorially as a weighted 
where $\gamma$ is an integer partition and $\mathfrak{d}$ is an \emph{arbitrarily labelled cylindric diagrams} (see Section \ref{alcd}).

In other words, Borodin's identity may be rewritten in the form:
\[
\sum_{\mathfrak{c} \in \CPP(\pi)} z^{|\mathfrak{c}|} =  \sum_{(\gamma,\mathfrak{d}) \in (\IP, \ALCD(\pi))} z^{|\mathfrak{d}| + T \, |\gamma|} 
\]

Here $\IP$ denotes the set of all integer partitions and $\ALCD(\pi)$ denotes the set of all arbitrarily labelled cylindric diagrams with profile $\pi$ (see Section \ref{alcd}).

Our goal is thus to find, for each possible profile $\pi$, a weight-preserving bijection between the
sets $\CPP(\pi)$ and the tuple $(\IP, \ALCD(\pi))$.

\begin{equation}\label{bijection-type} 
\psi_\pi : (\IP, \ALCD(\pi)) \to \CPP(\pi) 
\end{equation}

 Our bijection will be such that it actually proves the following refined identity:

\begin{multline} 
\sum_{\mathfrak{c} \in \CPP(\pi)} z_1^{|\mu_1|}z_2^{|\mu_2|} \cdots z_T^{|\mu_T|} = \\
 \sum_{(\gamma,\mathfrak{d}) \in (\IP, \ALCD(\pi))} z_1^{|\gamma| + |\diag(1)|_\mathfrak{d}} z_2^{|\gamma| + |\diag(2)|_\mathfrak{d}} \cdots z_T^{|\gamma| + |\diag(T)|_\mathfrak{d}} 
\end{multline}

The factorized form of the refined version of the identity may be obtained from equation \ref{borodin} via the following replacements:

\begin{align*} 
z^{nT} & \mapsto z_1^{n} z_2^{n} \cdots z_{T}^{n} \\
z^{j - i + nT} & \mapsto z_1^{n} z_2^n \cdots z_i^{n} z_{i+1}^{n+1} \cdots z_{j}^{n+1} z_{j+1}^n \cdots z_{T}^n \quad \quad \quad \quad \text{ when } i < j \\
z^{j-i +(n+1)T} & \mapsto z_1^{n+1} z_2^{n+1} \cdots z_j^{n+1} z_{j+1}^n + \cdots z_i^{n} z_{i+1}^{n+1} \cdots z_{T}^{n+1} \quad \text{ when } i > j
\end{align*}

\begin{definition}\label{strong-weight}
We shall say that a bijection $\psi_\pi : (\IP, \ALCD(\pi)) \to \CPP(\pi)$ is \emph{strongly weight preserving}
if whenever $\psi_\pi(\nu, \mathfrak{d}) = \mathfrak{c} = (\mu^0, \mu^1, \ldots, \mu^T) $ we have for all $1 \leq k \leq T$ that:
\[ |\mu^k| = |\gamma| + |\diag(k)|_\mathfrak{d} \]
\end{definition}

\begin{lemma}
Strongly weight preserving implies weight preserving.
\end{lemma}

\section{Symmetric Functions}
\label{sym-fun}

In this section we recall some of the theory of symmetric functions. Although the goal of this paper is to give a bijective proof of Borodin's identity, we shall see that the bijective proof is very closely related to the algebraic proof. They key idea in the algebraic proof are certain commutation relations between operators acting on symmetric functions. In the bijective proof, these commutation relations correspond to the local rules described in the next section.

Let $\Lambda$ denote the ring of symmetric functions over the field of rational numbers \cite{macdonald}.
Whenever possible we shall suppress in our notation any mention to the variables in which the functions are symmetric.
When we must mention the variables explicitly we shall make use of the \emph{plethystic notation} \cite{garsia}.

In the plethystic notation addition corresponds to the union of two sets and multiplication corresponds to the cartesian product. 
For example, we write:
\begin{equation} X = x_1 + x_2 + \cdots \end{equation}
to denote the set of variables $\{x_1, x_2, \ldots \}$. We also write:
\begin{equation} XY = (x_1 + x_2 + \cdots)(y_1 + y_2, \ldots) \end{equation}
to denote the set of variables $\{x_1 y_1, x_1 y_2, \ldots, x_2 y_1, \ldots x_2 y_2 \ldots \}$.

Let us denote the generating function for the complete symmetric functions by:
\begin{equation} 
 \Omega[Xz] = \prod_i \frac{1}{1-x_iz} = \sum_n h_n z^n 
\end{equation}

The \emph{Cauchy Kernel} is given by:
\begin{equation}
  \Omega[XY] = \prod_{i,j} \frac{1}{1-x_i y_j} = \sum_\lambda S_\lambda(X) S_\lambda(Y) 
\end{equation}

Recall that the \emph{Schur functions} are an orthonormal basis for $\Lambda$ with respect to the \emph{Hall Inner product}
\begin{equation} 
 \langle S_\lambda \,|\, S_\mu \rangle = \delta_{\lambda, \mu} 
\end{equation}

The operator $\Omega^*[Xz]$ is defined to be adjoint to the operator $\Omega[Xz]$ with respect to the Hall inner product.

\begin{equation*}
 \langle f(X) \, | \, \Omega^*[Xz] g(X) \rangle = \langle \Omega[Xz] f(X) \, | \, g(X) \rangle  
\end{equation*}

The Pieri formulae may be written in the form:

\begin{equation}\label{p1} \Omega[Xz] S_\mu[X] 
 = \sum_{\lambda \in U(\mu)} S_\lambda[X] z^{|\lambda|- |\mu|}\end{equation}
\begin{equation}\label{p2} \Omega^*[Xz]S_\lambda[X] = S_\lambda[X+z] = \sum_{\mu \in D(\mu)} S_\mu[X] z^{|\lambda|- |\mu|} \end{equation}

where $U(\mu)$ denotes the set of partitions which can be obtained from $\mu$ by adding a horizontal strip and $D(\lambda)$ denotes the set of partitions which can be obtained from $\lambda$ by removing a horizontal strip.

The following commutation relations are well known. They are essentially those of the Heisenberg algebra. They form the backbone of all algebraic proofs of hook length formulae for plane partitions:
\begin{lemma}
\boxed{ \Omega^*[Xu] \, \Omega[Xv] = \frac{1}{1-uv} \, \Omega[Xv] \,\Omega^*[Xu]}
\end{lemma}

We shall now sketch an algebraic proof of Borodin's identity.
We begin with a number of small lemmas. 
Let $D_z$ denote the `degree'' operator:
\begin{equation} \label{commutation}
D_z S_\lambda[X] = z^{|\lambda|} S_\lambda[X] 
\end{equation}

The degree operator satisfies the following commutation relations:
\begin{lemma} \label{degree}
\begin{align}
D_z \, \Omega[Xu] & = \Omega[Xuz] \, D_z \\
D_z \, \Omega^*[Xu] & = \Omega^*[Xuz^{-1}] \, D_z 
\end{align}
\begin{proof}
This fact follows immediately from the Pieri formulae.
\end{proof}
\end{lemma}

For notational convenience we shall define:
\begin{align}
G^0(z) & = \Omega[Xz] \\
G^1(z) & = \Omega^*[Xz]
\end{align}

\begin{lemma}  \label{lhs}
The left hand side of the refined version of equation \ref{maintheorem} may be expressed in the form:
\begin{equation} 
\lhs(\pi) = \sum_\mu \langle S_\mu \,|\, G^{\pi_0}(u_0) G^{\pi_1}(u_1)  \cdots G^{\pi_T}(u_T) D_w \, S_\mu  \rangle 
\end{equation}
where:
\begin{align} \label{specialization}
w & = z_0 z_1 \cdots z_{T-1} \\ 
u_k & = \begin{cases}
z_0 z_1 \cdots z_{k-1} & \text{ if $\pi_k = 1$} \\
z_0^{-1} z_1^{-1} \cdots z_{k-1}^{-1} & \text{ if $\pi_k = 0$}
\end{cases}\label{specialization2}
\end{align}

\begin{proof}
From the ``interlacing sequence'' definition of a cylindric plane partition it is clear that a cylindric plane partition is constructed
by successively adding and removing horizontal strips. 
The degree operator $D_z$ is used to keep track of the number of cubes in the resulting cylindric plane partition.

Using the fact that the Schur functions are orthonormal with respect to the hall inner product
we may write:
\begin{equation}
\lhs(\pi) = \sum_\mu \langle S_\mu \,|\, D_{z_0} \, G^{\pi_0}(1) \, D_{z_1} \, G^{\pi_1}(1)   \cdots D_{z_{T-1}} \, G^{\pi_T}(1) \, S_\mu  \rangle
\end{equation}

It remains to commute all the shift operators to the right hand side using Lemma \ref{degree}.
\end{proof}
\end{lemma}

Next let us define:
\begin{definition}\label{M-def}
\begin{equation} 
M_\pi(m)  = \sum_\mu \langle S_\mu \,|\, \prod_{\substack{k =1\\ \pi_k = 0}}^T \Omega [Xu_kw^m] 
\prod_{\substack{k =1 \\ \pi_k = 1}}^T \Omega^*[Xu_k] \, D_w \, S_\mu \rangle 
\end{equation}
\end{definition}

\begin{lemma} \label{cylindric-shift}
\[ M_\pi(m) =  \prod_{\substack{(i,j) \\ \pi_i \neq \pi_j}} \frac{1}{1-u_i u_j w^{m+1}} \, M_\pi(m+1) \]
\begin{proof}
This is a straightforward calculation. Using the fact that the Schur functions are orthogonal with respect to the Hall inner product, we may write:
\begin{align}
M_\pi(m) & = \sum_{\mu, \lambda} 
\langle S_\mu \,|\, \prod_{\substack{k =1\\ \pi_k = 0}}^T \Omega[Xu_kw^m] \, S_\lambda \rangle 
\langle S_\lambda \,|\,\prod_{\substack{k =1 \\ \pi_k = 1}}^T \Omega^*[Xu_k] \, D_w \, S_\mu \rangle \\
& = \sum_{\mu, \lambda} 
\langle S_\lambda \,|\,\prod_{\substack{k =1 \\ \pi_k = 1}}^T \Omega^*[Xu_k] \, D_w \, S_\mu \rangle 
\langle S_\mu \,|\, \prod_{\substack{k =1\\ \pi_k = 0}}^T \Omega[Xu_kw^m] \, S_\lambda \rangle \\
& = \sum_{\lambda} 
\langle S_\lambda \,|\,\prod_{\substack{k =1 \\ \pi_k = 1}}^T \Omega^*[Xu_k] \, 
D_w \,\prod_{\substack{k =1\\ \pi_k = 0}}^T \Omega[Xu_kw^m] \, P_\lambda \rangle_{q,t} 
\end{align}
Next applying the commutation relations of Lemma \ref{degree} and Lemma \ref{commutation} we have:
\begin{align}
M_\pi(m) & = \sum_{\lambda} 
\langle S_\lambda \,|\,\prod_{\substack{k =1 \\ \pi_k = 1}}^T \Omega^*[Xu_k] \, 
D_w \,\prod_{\substack{k =1\\ \pi_k = 0}}^T \Omega[Xu_kw^m] \, S_\lambda \rangle  \\ 
& = \sum_{\lambda} 
\langle S_\lambda \,|\,\prod_{\substack{k =1 \\ \pi_k = 1}}^T \Omega^*[Xu_k] \, 
\,\prod_{\substack{k =1\\ \pi_k = 0}}^T \Omega[Xu_kw^{m+1}] \, D_w \, S_\lambda \rangle \\
& = \prod_{\substack{(i,j) \\ \pi_i \neq \pi_j}} \frac{1}{1-u_i u_j w^{m+1}} \sum_{\lambda} 
\langle S_\lambda \,| \,\prod_{\substack{k =1\\ \pi_k = 0}}^T \Omega[Xu_kw^{m+1}] \, 
\prod_{\substack{k =1 \\ \pi_k = 1}}^T \Omega^*[Xu_k] \, D_w \, S_\lambda \rangle \\
& =  \prod_{\substack{(i,j) \\ \pi_i \neq \pi_j}} \frac{1}{1-u_i u_j w^{m+1}} \, M_\pi(m+1)
\end{align}
\end{proof}
\end{lemma}
In the limit we have:

\begin{lemma} \label{infinity}
\begin{equation} 
M_\pi(\infty)  = \prod_{n \geq 1} \frac{1}{1-w^n} 
\end{equation}
\begin{proof}
In order for this limit to even make sense, we must have $|z_i| < 1$ for all $i$, in which case:
\[ \lim_{m \to \infty} \Omega[X u_k \omega^m] = 1 \]

Since $\Omega^*[X u_k]$ is a degree lowering operator, it follows that:
\begin{align*}
 \lim_{m \to \infty} M_\pi(m) 
& =  \sum_\mu \langle S_\mu \,|\,  \prod_{\substack{k =1 \\ \pi_k = 1}}^T \Omega^*[Xu_k] \, D_w \, S_\mu \rangle \\
& = \sum_\mu \langle S_\mu | D_w \, S_\mu \rangle \\
& = \sum_\mu \omega^{|\mu|} \\
& = \prod_{n \geq 1} \frac{1}{1-w^n} 
\end{align*}

\end{proof}
\end{lemma}

The proof of the refined version of Theorem \ref{borodin} now proceeds as follows. We begin by applying Lemma \ref{lhs}

\begin{align*}
\phantom{=} \sum_{\mathfrak{c} \in \CPP(\pi)}  z^{|\mathfrak{p}|} 
& = \sum_\mu \langle S_\mu \,|\, G^{\pi_0}(u_0) G^{\pi_1}(u_1)  \cdots G^{\pi_T}(u_T) D_w \, S_\mu  \rangle \\
\end{align*}
Next we repeatedly applies the commutation relations of Lemma \ref{commutation}, followed by definition \ref{M-def}.
\begin{align*}
& = \prod_{\substack{i < j \\ \pi_i > \pi_j} }  \frac{1}{1-u_i u_j}
\sum_\mu \langle S_\mu \,|\, \prod_{\substack{k =1\\ \pi_k = 0}}^T \Omega[Xu_k] 
\prod_{\substack{k =1 \\ \pi_k = 1}}^T \Omega^*[Xu_k] \, D_w \, S_\mu \rangle \\
& = \prod_{\substack{i < j \\ \pi_i > \pi_j} } \frac{1}{1-u_i u_j} M_\pi(0) 
\end{align*}
We then repeatedly apply Lemma \ref{cylindric-shift}.
\begin{align*}
& = \prod_{\substack{i < j \\ \pi_i > \pi_j} } \frac{1}{1-u_i u_j} \prod_{m \geq 0} 
\left ( \prod_{\substack{(i,j) \\ \pi_i \neq \pi_j}} \frac{1}{1-u_i u_j w^{m+1}} \right ) M_\pi (\infty) 
\end{align*}
Splitting the second product into two, and combining it with the first we have:

\begin{align*}
& = \prod_{m \geq 1}
\left ( \prod_{\substack{i < j \\ \pi_i > \pi_j} } \frac{1}{1-u_i u_jw^{m-1}} \right )
\left ( \prod_{\substack{i > j \\ \pi_i > \pi_j} } \frac{1}{1-u_i u_jw^m}\right ) M_\pi (\infty)
\end{align*}
Finally, applying Lemma \ref{infinity} we have:
\begin{align*}
& = \prod_{m \geq 1} \frac{1}{1-w^m}
\left ( \prod_{\substack{i < j \\ \pi_i > \pi_j} } \frac{1}{1-u_i u_jw^{m-1}} \right )
\left ( \prod_{\substack{i > j \\ \pi_i > \pi_j} } \frac{1}{1-u_i u_jw^m}\right ) 
\end{align*}

To obtained the non-refined version of the Theorem, it suffices to take the following specialization of variables on both sides:
\begin{align} \label{specialization}
w & = z^{|T|} \\ 
u_k & = \begin{cases}
z^k & \text{ if $\pi_k = 1$} \\
z^{-k} & \text{ if $\pi_k = 0$}
\end{cases}\label{specialization2}
\end{align}

\section{Local Rule}\label{local-rule}

There is a lot that we won't say about the local rules because they have been discussed in detail elsewhere (see for example \cite{vanleeuwen}). Suffice is to say the local rules give a bijective proof of the commutation relation in lemma \ref{commutation} which we rewrite in the form:

\begin{equation}
\langle \, S_\alpha \, |  \,\Omega^*[Xu] \, \Omega[Xv] \, S_\beta  \,\rangle
= \frac{1}{1-uv} \, \langle  \,S_\alpha  \,|  \,\Omega[Xv] \,\Omega^*[Xu]  \,S_\beta  \,\rangle
\end{equation}

 That is, a local rule is a map with type signature:

\begin{equation} 
\mathfrak{D}_{\alpha,\beta} : U(\alpha) \cap U(\beta) \to (\mathbb{Z}_{\geq 0}, D(\alpha) \cap D(\beta))
\end{equation}

Such that if $(\ell, \nu) = \mathfrak{D}_{\alpha,\beta}(\lambda)$ then the following weight conditions are satisfied:

\begin{align} 
|\lambda / \alpha| & = |\beta / \nu| + \ell \\
|\lambda / \beta| & = |\alpha / \nu| + \ell 
\end{align}

These two weight conditions may be combined to give:
\begin{equation} 
|\lambda / \nu| = |\alpha / \nu| + |\beta / \nu| + \ell 
\end{equation}

Sometimes the local rule is represented graphically as follows:

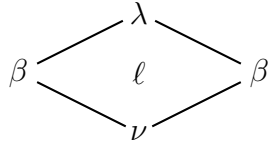
\begin{center}
\begin{tikzpicture}[scale=0.8]

\begin{scope}
\path (-2,1) node[shape=circle,inner sep=0.3mm](B0){$\beta$};
\path (0,0) node[shape=circle,inner sep=0.3mm](C0){$\nu$};
\path (2,1) node[shape=circle,inner sep=0.3mm](D0){$\beta$};
\path (0,2) node[shape=circle,inner sep=0.3mm](C1){$\lambda$};

\draw[thick]  (B0) -- (C0) -- (D0);

\draw[thick] (B0) -- (C1);
\draw[thick] (D0) -- (C1);

\path (0,1) node[shape=circle,inner sep=0.3mm](c0){$\ell$};

\end{scope}

\end{tikzpicture}
\end{center}

Observe that the weight condition implies immediately the following lemma:
\begin{lemma}\label{local-help2}
$\mathfrak{D}_{\gamma,\gamma}(\gamma) = (0,\gamma)$
\end{lemma}

The inverse local rule has type signature:

\begin{equation} 
\mathfrak{U}_{\alpha,\beta} : (\mathbb{Z}_{\geq 0}, D(\alpha) \cap D(\beta)) \to U(\alpha) \cap U(\beta)  
\end{equation}

Fomin's original local rule \cite{fomin1} corresponded to the RS correspondence between permutations and pairs of standard tableaux rather than the full RSK correspondence between integer matrices and pairs of semi-standard tableaux. 

There are, in fact, two equally natural versions of Fomin's original local rule which differ only by conjugation. As a consequence there are two equally naturally versions of the generalized local rule, one of which gives the RSK correspondence and other which gives the Burge correspondence \cite{burge}. These two local rules are related by the Schutzenburger involution \cite{vanleeuwen}.

We shall describe here only the local rule associated to the Burge correspondence, because it is the rule used in our example in section \ref{growth}. The reader is refered to \cite{vanleeuwen} for more details.

In the case of the Berge correspondence, the operator $\mathfrak{D}_{\alpha,\beta}$ is defined as follows.
Suppose that $(\ell,\nu) = \mathfrak{D}_{\alpha,\beta}(\lambda)$. Let $\overline{A}$ denote the set of columns of $\lambda$ which are longer than the corresponding columns of $\alpha$
and let $\overline{B}$ denote the set of columns of $\lambda$ which are longer than the corresponding columns of $\beta$.

Next, for each $i \in \overline{A} \cap \overline{B} $ let $\delta(i) \not\in \overline{A} \cup \overline{B}$ denote the largest integer such that $\delta(i) < i$
and $\delta(i) \neq \delta(j)$ for any $j \in \overline{A} \cap \overline{B}$ such that $j > i$.
Let:
\[ \overline{C} = \{\delta(i) > 0 \, | \, i \in \overline{A} \cap \overline{B} \} \]
Finally let $\nu$ be the partition obtained from $\lambda$ by removing a box from the end of each of the columns indexed by $\overline{A} \cup \overline{B} \cup \overline{C}$ and let
\[\ell = \#\{ \delta(i) \leq 0 \, | \, i \in \overline{A} \cap \overline{B} \} \]

Here is an example:
\[ \boxed{\mathfrak{D}_{(6,5,5,3),(6,6,5,2)}(7,6,5,3,1) = (1,(6,5,4,2))}\]

The calculation proceeds as follows:

\begin{align*}
\lambda' & = (5,4,4,3,3,2,1) \\
\overline{A} = \cols(\lambda / \alpha) & = \{1,6,7\} \\
\overline{B} = \cols(\lambda / \beta) & = \{1,3,7\} \\
\end{align*}

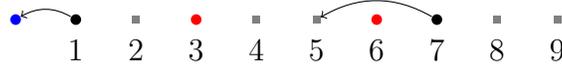
\tikzstyle{blacknode}=[shape=circle,fill,inner sep=0.5mm]
\tikzstyle{rednode}=[shape=circle,red, fill, inner sep=0.5mm]
\tikzstyle{graynode}=[shape=rectangle, gray, fill, inner sep=0.5mm]
\tikzstyle{bluenode}=[shape=circle, blue, fill, inner sep=0.5mm]

\begin{center}\begin{tikzpicture}[scale=0.8]

\path (-1,0) node[bluenode](n0){};

\path (0,0) node[blacknode](n1){};
\path (1,0) node[graynode](n2){};
\path (2,0) node[rednode](n3){};
\path (3,0) node[graynode](n4){};
\path (4,0) node[graynode](n5){};
\path (5,0) node[rednode](n6){};
\path (6,0) node[blacknode](n7){};
\path (7,0) node[graynode](n8){};
\path (8,0) node[graynode](n9){};

\path (0,-0.5) node[](){1};
\path (1,-0.5) node[](){2};
\path (2,-0.5) node[](){3};
\path (3,-0.5) node[](){4};
\path (4,-0.5) node[](){5};
\path (5,-0.5) node[](){6};
\path (6,-0.5) node[](){7};
\path (7,-0.5) node[](){8};
\path (8,-0.5) node[](){9};

\path[<-] (n5) edge [bend left] (n7);
\path[<-] (n0) edge [bend left] (n1);
\end{tikzpicture}\end{center}

\begin{align*}
\overline{C} & = \{5\} \\
\mu' & = (4,4,3,3,2,1) \\
\ell & = 1 
\end{align*}

The inverse operator $\mathfrak{U}_{\alpha,\beta}$ is defined similarly.
Suppose that $\lambda = \mathfrak{U}_{\alpha,\beta}(\ell,\nu)$. Let $A$ denote the set of columns of $\alpha$ which are longer than the corresponding columns of $\nu$
and let $B$ denote the set of columns of $\beta$ which are longer than the corresponding columns of $\nu$.

Next, for each $i \in A \cap B$ let $\epsilon(i) \not \in A \cup B$ be the smallest integer such that $\epsilon(i) > i$ and $\epsilon(i) \neq \epsilon(j)$ for any $j \in A \cap B$ with $j < i$.
Let 
\[ C = \{ \epsilon(i) \,|\, i \in A \cap B \}\] 
Finally let $D$ denote the first $\ell$ elements of the complement of the set $A \cup B \cup C$ 
and let $\lambda$ be the partition obtained from $\nu$ by adding a box to the end of each of the columns in $A \cup B \cup C \cup D$.

Here is the inverse of our example:

\[ \boxed{\mathfrak{U}_{(6,5,5,3),(6,6,5,2)}(1,(6,5,4,2)) = (7,6,5,3,1) }\]

The calculation is straightforeward:

\begin{align*}
\nu' & = (4,4,3,3,2,1) \\
A = \cols(\alpha / \nu) & = \{3,5\} \\
B = \cols(\beta / \nu) & = \{5, 6\} \\
m & = 1
\end{align*}

\tikzstyle{blacknode}=[shape=circle,fill,inner sep=0.5mm]
\tikzstyle{rednode}=[shape=circle,red, fill, inner sep=0.5mm]
\tikzstyle{graynode}=[shape=rectangle, gray, fill, inner sep=0.5mm]
\tikzstyle{bluenode}=[shape=circle, blue, fill, inner sep=0.5mm]

\begin{center}\begin{tikzpicture}[scale=0.8]

\path (-1,0) node[bluenode](n0){};
\path (0,0) node[graynode](n1){};
\path (1,0) node[graynode](n2){};
\path (2,0) node[rednode](n3){};
\path (3,0) node[graynode](n4){};
\path (4,0) node[blacknode](n5){};
\path (5,0) node[rednode](n6){};
\path (6,0) node[graynode](n7){};
\path (7,0) node[graynode](n8){};
\path (8,0) node[graynode](n9){};

\path (0,-0.5) node[](){1};
\path (1,-0.5) node[](){2};
\path (2,-0.5) node[](){3};
\path (3,-0.5) node[](){4};
\path (4,-0.5) node[](){5};
\path (5,-0.5) node[](){6};
\path (6,-0.5) node[](){7};
\path (7,-0.5) node[](){8};
\path (8,-0.5) node[](){9};

\path[->] (n5) edge [bend left] (n7);
\path[->] (n0) edge [bend left] (n1);
\end{tikzpicture}\end{center}

\begin{align*}
C & = \{7\} \\
\lambda' & = (5,4,4,3,3,2,1) \\ 
\end{align*}

\section{Local Rule as Higher Order Function}\label{higher-order}

\subsection{Type Signatures}

Up until now we have defined the local rule as a simple map with type signature:
\[\mathfrak{U}_{\alpha,\beta} : U(\alpha) \cap U(\beta) \to (\mathbb{Z}_{\geq 0}, D(\alpha) \cap D(\beta)) \]
whose inverse map has type signature:
\[\mathfrak{D}_{\alpha,\beta} : (\mathbb{Z}_{\geq 0}, D(\alpha) \cap D(\beta)) \to U(\alpha) \cap U(\beta) \]

At slight risk of confusion, we shall also use the term ``local rule'' to refer to a certain \emph{higher order function}, in the sense of functional programming \cite{haskell}.
  
In functional programming, a higher order function is a function which takes as input a function, and returns as output a different function. 

Recall that our goal is to construct a weight preserving bijection, for each binary string $\pi$ between, on the one hand, the set of cylindric plane partitions $\CPP(\pi)$ with profile $\pi$, and on the other hand the pair $(\IP, \ALCD(\pi))$ where $\IP$ denotes the set of all integer partitions, and $\ALCD(\pi)$ denotes the set of arbitrarily labelled cylindric diagrams (see Section \ref{statement}).

For any $\pi \prec \pi'$ such that $\pi_i \neq \pi'_i$, the input function for our local rule $\mathfrak{L}_i$ will be a weight-preserving bijection of the form:
\[ \psi_\pi : (\IP, \ALCD(\pi)) \to \CPP(\pi) \]
while the output function is a weight preserving bijection of the form:
\[ \psi_{\pi'} : (\IP, \ALCD(\pi')) \to \CPP(\pi') \]
That is to say, the local rule $\mathfrak{L}_i$ will have type signature:

\[ \mathfrak{L}_i : ((\IP, \ALCD(\pi)) \to \CPP(\pi)) \to ((\IP, \ALCD(\pi')) \to \CPP(\pi')) \]
In other words:
\[ \mathfrak{L_i}[ \psi_\pi ] = \psi_{\pi'} \]

Let $\varphi_\pi$ and $\varphi_{\pi'}$ denote the inverse of $\psi_\pi$ and $\psi_{\pi'}$ respectively.
The ``inverse local rule'' $\mathfrak{M}_i$ is the higher order function with type signature:
\[ \mathfrak{M}_i : (\CPP(\pi) \to (\IP, \ALCD(\pi)) ) \to (\CPP(\pi') \to (\IP, \ALCD(\pi'))) \]
That is:
\[ \mathfrak{M}_i[ \varphi_{\pi} ] = \varphi_{\pi'} \]
Note that $\mathfrak{M}_i$ is only ``inverse'' to $\mathfrak{L}_i$ in the sense that:
\begin{align*}
\mathfrak{L}_i [ \psi_\pi ] \circ \mathfrak{M}_i [ \varphi_{\pi} ] & = \mathfrak{1}_{\CPP(\pi')} \\
\mathfrak{M}_i [ \varphi_{\pi} ] \circ \mathfrak{L}_i [ \psi_\pi ] & = \mathfrak{1}_{(\IP, \ALCD(\pi'))} 
\end{align*}
It is not possible to compose $\mathfrak{L}_i$ and $\mathfrak{M}_i$ directly due to incompatible type signatures.

\subsubsection{Adding and removing boxes}

An \emph{inside corner box} of an arbitrarily labelled cylindric plane diagram $\mathfrak{d}'$ with profile $\pi'$ is a box with cylindric inversion coordinates $(i,i+1,0)$ where $i$ is an inversion in the profile $\pi'$. 

An inside corner box of an arbitrarily labelled cylindric diagram $\mathfrak{d}'$ with profile $\pi'$ can always be removed to obtain an arbitrarily labelled cylindric diagram with profile $\pi$ where $\pi \prec \pi'$ (see Section \ref{profiles}). 

We shall denote this operator by:
\[ \mathfrak{l}_i : \ALCD(\pi') \to (\mathbb{Z}_{\geq 0}, \ALCD(\pi)) \]

Conversely, if $\mathfrak{d}$ is an arbitrarily labelled cylindric diagram with profile $\pi$ such that $\pi_i = 0$ and $\pi_i = 1$, then given an integer $m$ we may create a new arbitrarily labelled cylindric diagram $\mathfrak{d}'$ with profile $\pi' \succ \pi$ by adding a box with cylindric inversion coordinates $(i,i+1,0)$ and label $m$.

We shall denote this operator by:
\[ \mathfrak{r}_i : (\mathbb{Z}_{\geq 0}, \ALCD(\pi)) \to \ALCD(\pi') \]

\subsection{Definition of Local Rule}\label{local-def}

Choose any $(\gamma, \mathfrak{d}') \in (\IP, \ALCD(\pi'))$ and let 
\[ (m, \mathfrak{d}) = \mathfrak{l}[\mathfrak{d}'] \]

Suppose that:
\[ \psi_\pi (\gamma, \mathfrak{d}) = \mathfrak{c} = (\mu^0, \mu^1, \ldots, \mu^T) \]
Let: 
\begin{align*}
\alpha & = \mu^{i-1} \\
\gamma & = \mu^i \\
\beta & = \mu^{i+1} \\
\end{align*}
and let 
\[ \boxed{\lambda = \mathfrak{U}_{\alpha,\beta} (\gamma, m)} \]
We define:
\[ \mathfrak{L}_i[\psi_\pi] (\gamma, \mathfrak{d}') = \mathfrak{c}' = (\mu^0, \ldots, \mu^{i-1}, \lambda, \mu^{i+1}, \ldots \mu^T) \] 

Note that the horizontal strip condition in the definition of $\mathfrak{U}_{\alpha,\beta}$ ensures that this definition is well-defined.

The inverse local rule is defined similarly. Choose any cylindrical plane partition $\mathfrak{c}' = (\mu^0, \mu^1, \ldots, \mu^T)$ with profile $\pi'$.
Let us define:
\begin{align*}
\alpha & = \mu^{i-1} \\
\lambda & = \mu^i \\
\beta & = \mu^{i+1} \\
\end{align*}
Next let 
\[ \boxed{(m,\nu) = \mathfrak{D}_{\alpha,\beta}(\lambda)} \]
and let $\mathfrak{c}$ be the cylindric plane partition with profile $\pi$
given by $\mathfrak{c} = (\mu^0, \ldots, \mu^{i-1}, \nu, \mu^{i+1}, \ldots \mu^T)$.

If $\varphi_\pi(\mathfrak{c}) = (\gamma, \mathfrak{d})$ then we define 
\[ \mathfrak{M}_i[\varphi_\pi](\mathfrak{c}') = (\gamma, \mathfrak{d}')\]
where $\mathfrak{d}' = \mathfrak{r}_i[\mathfrak{d}]$

\newpage

\section{The Bijection}

In this section we construct recursively, for each possible profile $\pi$, the strongly weight-preserving bijection defined in equation \ref{bijection-type}: 

\[ \psi_\pi : (\IP, \ALCD(\pi)) \to \CPP(\pi) \]

\subsection{Idea of bijection}

In the special case when the arbitrarily labelled cylindric diagram $\mathfrak{d}$ has depth zero (see definition \ref{depth}), the strongly weight-preserving bijection is particularly simple:
\[ \psi_{\pi}(\gamma, \emptyset) = (\gamma,\gamma, \ldots \gamma) \]

The idea is to recurseively construct bijections, starting from this base case, by repeated application of the local rule (section \ref{higher-order}) and the rotation operator (definition \ref{shift}). 

Each application of the local rule corresponds to a commutation relation in the algebraic proof of Borodin's identity while each application of the rotation operator corresponds to the operation of splitting the trace operator into two and then interchanging the sum (see section \ref{sym-fun}) 

The recursion is not only over the number of inversions in the profile, but also over the \emph{depth} of the arbitrarily labelled cylindric diagram upon which the bijection is acting. Although the precise structure of the recursion is a little complicated to describe in words, it may be neatly encoded in a geometric object called a \emph{cylindric growth diagram}.

The local rules alone would not allow us to define the bijection on cylindric diagrams with depth greater than or equal than two. For this we must use the rotation operator.

\subsection{Cylindric growth diagrams}\label{growth}

The idea of a \emph{growth diagram} was first introduced by Fomin \cite{fomin1, fomin2}. Krattenthaler \cite{krattenthaler} made use of this framework to give a new bijective proof of Stanley's identity (equation \ref{stanley}). In the cylindric case we change the underlying poset, but the essential idea remains the same.

\begin{definition}
The \emph{cylindric poset} $\mathfrak{G}(n,m)$ is the quotient of $\mathbb{Z}^2$ via the equivalence relation:
\begin{equation}
(x,y) \equiv (x+n,y-m)
\end{equation}
\end{definition}

We shall write $v \lhd w$ to indicate that the vertex $v$ covers the vertex $w$ in the cylindric poset.

\begin{definition}
For any binary string $\pi$ containing $n$ zeros and $m$ ones, a \emph{path} in a cylindric growth diagram $\mathfrak{G}(n,m)$ with profile $\pi$ is a sequence of vertices: 
\[ p = (v_0, v_1, \ldots, v_{n+m-1}, v_{n+m}) \qquad \text{ with } \qquad v_0 = v_{n+m}\] 
 satisfying $v_{k-1} \lhd v_k$ if $\pi_k = 0$, otherwise $v_{k-1} \rhd v_k$.
\end{definition}

The partial order on $\bin(n,m)$ described in section \ref{profiles} induces a partial order in the set of paths of $\mathfrak{G}(n,m)$. The covering relations are given by $p \prec q$ if $\pi(p) \prec \pi(q)$ and there is only one vertex in the path $p$ which is not also in the path $q$. Here $\pi(p)$ denotes the profile of the path $p$ while $\pi(q)$ denotes the profile of the path $q$.

\begin{definition}
A \emph{face} in the cylindric poset is a set of four vertices $(u ; v1, v2 ; w) $ satisfying
$u \lhd v1 \lhd w$ and $u \lhd v2 \lhd w$. 
\end{definition}

We say that the face $(u ; v1, v2; w)$ lies \emph{above} the vertex $u$ and \emph{below} the vertex $w$

\begin{lemma}
There is a natural bijection between the cylindric diagram with profile $\pi \in \bin(n,m)$, and the subset of the cylindric poset $\mathfrak{G}(n,m)$ which lies below any given path with profile $\pi$.
\end{lemma}

This bijection maps the \emph{boxes} of the cylindric diagram (section \ref{diagrams}) to the \emph{faces} of the cylindric poset.


\begin{definition}\label{growth-def}
A \emph{cylindric growth diagram} with profile $\pi \in \bin(n,m)$ is a subset of cylindric diagram $\mathfrak{G}(n,m)$ which lies below a path with profile $\pi$ (refered to as the upper boundary), whose vertices are labelled by integer partitions, and whose faces are labelled by non-negative integers, in such that the following three conditions are satisfied:

\begin{enumerate}
\item If $v \lhd w$, and if $\lambda$ is the integer partition labelling the vertex $v$ and $\mu$ is the integer partition labelling the vertex $w$, then $\lambda / \mu$ is a horizontal strip.
\item All but finitely many vertices are labelled with the same integer partition $\gamma$.
\item If $( u; v1, v2; w)$ is a face with label $\ell$, and if the labels of $u$, $v1$, $v2$ and $w$ are $\mu$, $\alpha$, $\beta$ and $\lambda$ respectively then:
\[ \mathfrak{U}_{\alpha,\beta}(\ell,\mu) = \lambda \]
\end{enumerate}
\end{definition}

The invertibility of the local rule implies that condition $(3)$ could be equivalently formulated as follows:
\begin{lemma}
If $( u; v1, v2; w)$ is a face of a cylindric diagram $\mathfrak{G}(n,m)$ with label $\ell$, and if the labels of $u$, $v1$, $v2$ and $w$ are $\mu$, $\alpha$, $\beta$ and $\lambda$ respectively then:
\[ \mathfrak{D}_{\alpha,\beta}(\lambda) = (\ell,\mu) \]
\end{lemma}

Here is an example:

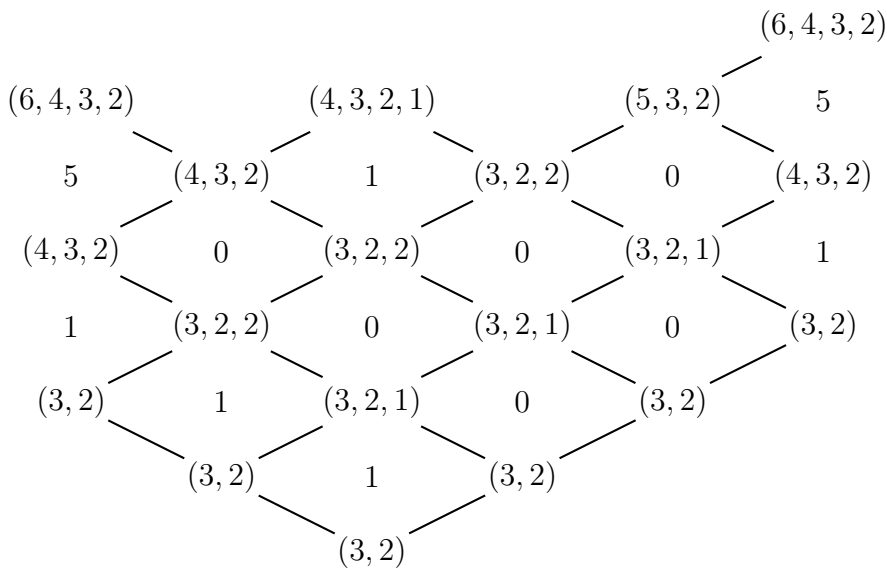
\begin{center}
\begin{tikzpicture}[scale=1]

\begin{scope}
\path (-4,2) node[shape=circle,inner sep=0.3mm](A0){$(3,2)$};
\path (-2,1) node[shape=circle,inner sep=0.3mm](B0){$(3,2)$};
\path (0,0) node[shape=circle,inner sep=0.3mm](C0){$(3,2)$};
\path (2,1) node[shape=circle,inner sep=0.3mm](D0){$(3,2)$};
\path (4,2) node[shape=circle,inner sep=0.3mm](E0){$(3,2)$};
\path (6,3) node[shape=circle,inner sep=0.3mm](F0){$(3,2)$};

\path (-4,4) node[shape=circle,inner sep=0.3mm](A1){$(4,3,2)$};
\path (-2,3) node[shape=circle,inner sep=0.3mm](B1){$(3,2,2)$};
\path (0,2) node[shape=circle,inner sep=0.3mm](C1){$(3,2,1)$};
\path (2,3) node[shape=circle,inner sep=0.3mm](D1){$(3,2,1)$};
\path (4,4) node[shape=circle,inner sep=0.3mm](E1){$(3,2,1)$};
\path (6,5) node[shape=circle,inner sep=0.3mm](F1){$(4,3,2)$};

\path (-4,6) node[shape=circle,inner sep=0.3mm](A2){$(6,4,3,2)$};
\path (-2,5) node[shape=circle,inner sep=0.3mm](B2){$(4,3,2)$};
\path (0,4) node[shape=circle,inner sep=0.3mm](C2){$(3,2,2)$};
\path (2,5) node[shape=circle,inner sep=0.3mm](D2){$(3,2,2)$};
\path (4,6) node[shape=circle,inner sep=0.3mm](E2){$(5,3,2)$};
\path (6,7) node[shape=circle,inner sep=0.3mm](F2){$(6,4,3,2)$};

\path (0,6) node[shape=circle,inner sep=0.3mm](C3){$(4,3,2,1)$};

\draw[thick] (A0) -- (B0) -- (C0) -- (D0) -- (E0) -- (F0);
\draw[thick] (A1) -- (B1) -- (C1) -- (D1) -- (E1) -- (F1);
\draw[thick] (A2) -- (B2) -- (C2) -- (D2) -- (E2) -- (F2);

\draw[thick] (A1) -- (B2) -- (C3);
\draw[thick] (A0) -- (B1) -- (C2);
\draw[thick] (B0) -- (C1);
\draw[thick] (D0) -- (C1);
\draw[thick] (E0) -- (D1) -- (C2);
\draw[thick] (F0) -- (E1) -- (D2) -- (C3);
\draw[thick] (F1) -- (E2);

\path (-4,3) node[shape=circle,inner sep=0.3mm](a0){$1$};
\path (-2,2) node[shape=circle,inner sep=0.3mm](b0){$1$};
\path (0,1) node[shape=circle,inner sep=0.3mm](c0){$1$};
\path (2,2) node[shape=circle,inner sep=0.3mm](dD0){$0$};
\path (4,3) node[shape=circle,inner sep=0.3mm](e0){$0$};
\path (6,4) node[shape=circle,inner sep=0.3mm](f0){$1$};

\path (-4,5) node[shape=circle,inner sep=0.3mm](a1){$5$};
\path (-2,4) node[shape=circle,inner sep=0.3mm](b1){$0$};
\path (0,3) node[shape=circle,inner sep=0.3mm](c1){$0$};
\path (2,4) node[shape=circle,inner sep=0.3mm](d1){$0$};
\path (4,5) node[shape=circle,inner sep=0.3mm](e1){$0$};
\path (6,6) node[shape=circle,inner sep=0.3mm](f1){$5$};

\path (0,5) node[shape=circle,inner sep=0.3mm](v2){$1$};
\end{scope}

\end{tikzpicture}
\end{center}

Note that we have truncated the diagram below the lower boundary, where all vertices have the same label.

\begin{lemma}
For any path $p$ in a cylindric growth diagram $\mathfrak{G}(n,m)$ with profile $\pi$, the sequence of profiles associated
to the vertices of the path form a cylindric plane partition.
\end{lemma}

In particular, the sequence of partitions labelling the vertices along the upper boundary form an element of $\CPP(\pi)$.

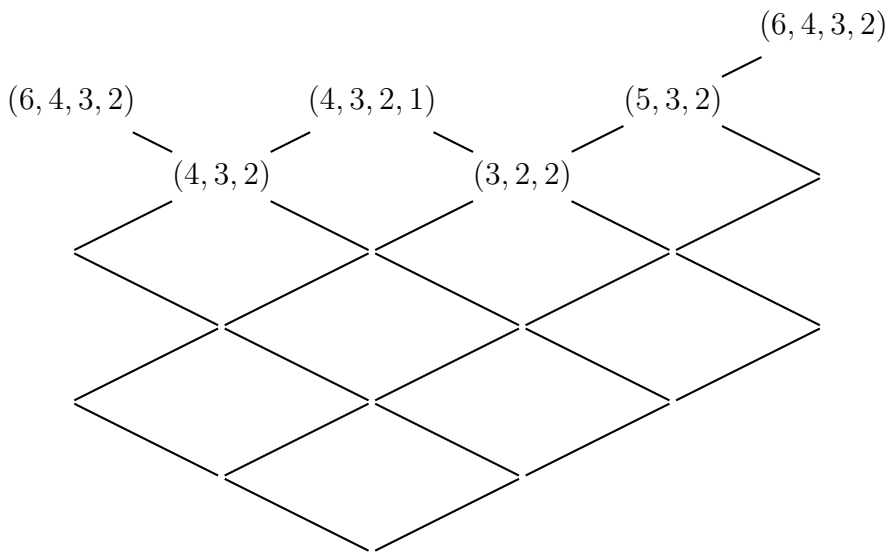
\begin{center}
\begin{tikzpicture}[scale=1]

\begin{scope}
\path (-4,2) node[shape=circle,inner sep=0.3mm](A0){};
\path (-2,1) node[shape=circle,inner sep=0.3mm](B0){};
\path (0,0) node[shape=circle,inner sep=0.3mm](C0){};
\path (2,1) node[shape=circle,inner sep=0.3mm](D0){};
\path (4,2) node[shape=circle,inner sep=0.3mm](E0){};
\path (6,3) node[shape=circle,inner sep=0.3mm](F0){};

\path (-4,4) node[shape=circle,inner sep=0.3mm](A1){};
\path (-2,3) node[shape=circle,inner sep=0.3mm](B1){};
\path (0,2) node[shape=circle,inner sep=0.3mm](C1){};
\path (2,3) node[shape=circle,inner sep=0.3mm](D1){};
\path (4,4) node[shape=circle,inner sep=0.3mm](E1){};
\path (6,5) node[shape=circle,inner sep=0.3mm](F1){};

\path (-4,6) node[shape=circle,inner sep=0.3mm](A2){$(6,4,3,2)$};
\path (-2,5) node[shape=circle,inner sep=0.3mm](B2){$(4,3,2)$};
\path (0,4) node[shape=circle,inner sep=0.3mm](C2){};
\path (2,5) node[shape=circle,inner sep=0.3mm](D2){$(3,2,2)$};
\path (4,6) node[shape=circle,inner sep=0.3mm](E2){$(5,3,2)$};
\path (6,7) node[shape=circle,inner sep=0.3mm](F2){$(6,4,3,2)$};

\path (0,6) node[shape=circle,inner sep=0.3mm](C3){$(4,3,2,1)$};

\draw[thick] (A0) -- (B0) -- (C0) -- (D0) -- (E0) -- (F0);
\draw[thick] (A1) -- (B1) -- (C1) -- (D1) -- (E1) -- (F1);
\draw[thick] (A2) -- (B2) -- (C2) -- (D2) -- (E2) -- (F2);

\draw[thick] (A1) -- (B2) -- (C3);
\draw[thick] (A0) -- (B1) -- (C2);
\draw[thick] (B0) -- (C1);
\draw[thick] (D0) -- (C1);
\draw[thick] (E0) -- (D1) -- (C2);
\draw[thick] (F0) -- (E1) -- (D2) -- (C3);
\draw[thick] (F1) -- (E2);

\path (-4,3) node[shape=circle,inner sep=0.3mm](a0){};
\path (-2,2) node[shape=circle,inner sep=0.3mm](b0){};
\path (0,1) node[shape=circle,inner sep=0.3mm](c0){};
\path (2,2) node[shape=circle,inner sep=0.3mm](dD0){};
\path (4,3) node[shape=circle,inner sep=0.3mm](e0){};
\path (6,4) node[shape=circle,inner sep=0.3mm](f0){};

\path (-4,5) node[shape=circle,inner sep=0.3mm](a1){};
\path (-2,4) node[shape=circle,inner sep=0.3mm](b1){};
\path (0,3) node[shape=circle,inner sep=0.3mm](c1){};
\path (2,4) node[shape=circle,inner sep=0.3mm](d1){};
\path (4,5) node[shape=circle,inner sep=0.3mm](e1){};
\path (6,6) node[shape=circle,inner sep=0.3mm](f1){};

\path (0,5) node[shape=circle,inner sep=0.3mm](v2){};
\end{scope}

\end{tikzpicture}
\end{center}

Note that this is just a rotation of the example cylindric plane partition given in setion \ref{cpp}.

\begin{lemma}
If the labels of the vertices of a cylindric growth diagram $\mathfrak{G}(n,m)$ with profile $\pi$ are forgotten, then we obtain an arbitrarily labelled cylindric diagram. 
\end{lemma}

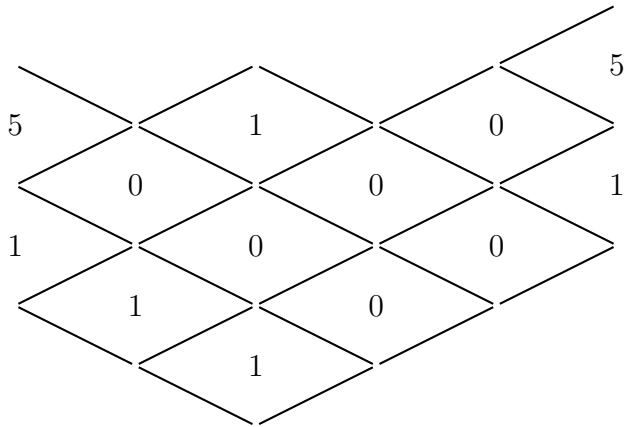
\begin{center}
\begin{tikzpicture}[scale=0.8]

\begin{scope}
\path (-4,2) node[shape=circle,inner sep=0.3mm](A0){};
\path (-2,1) node[shape=circle,inner sep=0.3mm](B0){};
\path (0,0) node[shape=circle,inner sep=0.3mm](C0){};
\path (2,1) node[shape=circle,inner sep=0.3mm](D0){};
\path (4,2) node[shape=circle,inner sep=0.3mm](E0){};
\path (6,3) node[shape=circle,inner sep=0.3mm](F0){};

\path (-4,4) node[shape=circle,inner sep=0.3mm](A1){};
\path (-2,3) node[shape=circle,inner sep=0.3mm](B1){};
\path (0,2) node[shape=circle,inner sep=0.3mm](C1){};
\path (2,3) node[shape=circle,inner sep=0.3mm](D1){};
\path (4,4) node[shape=circle,inner sep=0.3mm](E1){};
\path (6,5) node[shape=circle,inner sep=0.3mm](F1){};

\path (-4,6) node[shape=circle,inner sep=0.3mm](A2){};
\path (-2,5) node[shape=circle,inner sep=0.3mm](B2){};
\path (0,4) node[shape=circle,inner sep=0.3mm](C2){};
\path (2,5) node[shape=circle,inner sep=0.3mm](D2){};
\path (4,6) node[shape=circle,inner sep=0.3mm](E2){};
\path (6,7) node[shape=circle,inner sep=0.3mm](F2){};

\path (0,6) node[shape=circle,inner sep=0.3mm](C3){};

\draw[thick] (A0) -- (B0) -- (C0) -- (D0) -- (E0) -- (F0);
\draw[thick] (A1) -- (B1) -- (C1) -- (D1) -- (E1) -- (F1);
\draw[thick] (A2) -- (B2) -- (C2) -- (D2) -- (E2) -- (F2);

\draw[thick] (A1) -- (B2) -- (C3);
\draw[thick] (A0) -- (B1) -- (C2);
\draw[thick] (B0) -- (C1);
\draw[thick] (D0) -- (C1);
\draw[thick] (E0) -- (D1) -- (C2);
\draw[thick] (F0) -- (E1) -- (D2) -- (C3);
\draw[thick] (F1) -- (E2);

\path (-4,3) node[shape=circle,inner sep=0.3mm](a0){$1$};
\path (-2,2) node[shape=circle,inner sep=0.3mm](b0){$1$};
\path (0,1) node[shape=circle,inner sep=0.3mm](c0){$1$};
\path (2,2) node[shape=circle,inner sep=0.3mm](dD0){$0$};
\path (4,3) node[shape=circle,inner sep=0.3mm](e0){$0$};
\path (6,4) node[shape=circle,inner sep=0.3mm](f0){$1$};

\path (-4,5) node[shape=circle,inner sep=0.3mm](a1){$5$};
\path (-2,4) node[shape=circle,inner sep=0.3mm](b1){$0$};
\path (0,3) node[shape=circle,inner sep=0.3mm](c1){$0$};
\path (2,4) node[shape=circle,inner sep=0.3mm](d1){$0$};
\path (4,5) node[shape=circle,inner sep=0.3mm](e1){$0$};
\path (6,6) node[shape=circle,inner sep=0.3mm](f1){$5$};

\path (0,5) node[shape=circle,inner sep=0.3mm](v2){$1$};
\end{scope}

\end{tikzpicture}
\end{center}

Note that the arbitrarily labelled cylindric diagram is a rotation of the example given in section \ref{alcd}.

A growth diagram should be thought of as an object which interpolates between the LHS and the RHS of the bijection which we wish to establish.

\begin{proposition}
To every cylindric plane partition, there is a uniquely associated cylindric growth diagram.
\begin{proof}
Once the labels on the upper boundary have been specified, property (3) of definition \ref{growth-def} ensures that there is a unique way in which to label the remaining faces and vertices.
\end{proof}
\end{proposition}

\begin{proposition}
To every pair $(\gamma, \mathfrak{d})$ where $\gamma$ is an integer partition and $\mathfrak{d}$ is an arbitrarily labelled cylindric diagram, there exists a uniquely defined cylindric growth diagram.
\begin{proof}
Let $d$ denote the depth of $\mathfrak{d}$. To every face with depth greater than $d$ assign the label $0$. To every vertex lying below a face with depth greater than $d$, assign the label $\gamma$. Property $(3)$ of definition \ref{growth-def} ensures that there is a unique way to label the remaining faces and vertices. Lemma \ref{local-help2} ensures that the resulting growth diagram is well-defined.
\end{proof}
\end{proposition}

\subsection{Higher order functions again}

Every path in a cylindric growth diagram corresponds to a cylindric plane partition. The ``higher order'' local rules $\mathfrak{M}_i$ allow us to lift a bijection from the cylindric plane partition associated to the path $p$ to a to the cylindric plane partition associated to the path $q$ when $q \succ p$ (see Section \ref{higher-order}). That is to say, each face in the cylindric growth diagram corresponds to a case of the ''simple'' local rule as described in section \ref{local-rule}.

The following lemma guarantees that, when there are multiple inversions in the profile string, the order in which the local rules are applied is of no importance.

\begin{lemma}\label{localcommute}
If $\pi$ has inversions at both positions $i$ and $j$ then:
\[ \mathfrak{L}_i \circ \mathfrak{L}_j [ \psi_\pi ] = \mathfrak{L}_j \circ \mathfrak{L}_i [ \psi_\pi ] \]
\begin{proof}
Without loss of generality we may assume that $j > i$. If $\pi$ has inversions at both positions $i$ and $j$ then $\pi_i = 0 = \pi_j$ and $\pi_{i+1} = 1 = \pi_{j+1}$
thus $j-i \geq 2$. Application of the local rule $\mathfrak{L}_i$ does not effect the $(j-1)$th diagonal. Similarly, application
of the local rule $\mathfrak{L}_j$ does not effect the $(i+1)$th diagonal, thus the two operators commute. 
\end{proof}
\end{lemma}

Given a path $p$ with profile $\pi_{\max}$ it is not possible to apply the local rule $\mathfrak{M}_i$ for any $i$. It is however possible to rotate the cylinder, and thus obtain a new path $\sigma(p)$ with profile $\pi_{\min}$ (see Section \ref{shift}).

If $\mathfrak{d}$ is an arbitrarily labelled cylindric diagram with profile $\pi$ and depth $d > 0$ then in order to construct the bijection $\psi_\pi(\gamma, \mathfrak{d})$ the cylindric shift operator will have to be applied $d-1$ times.

\section{The Weight}

In this section we prove that our bijection is strongly weight preserving (see definition \ref{strong-weight}).
The reader is advised to review the definitions in section \ref{local-def} before proceeding with the proof.

\begin{proposition}
If $(m, \mathfrak{d}) = \mathfrak{l}_i[\mathfrak{d}']$ then:
\[ |\diag(i)|_{\mathfrak{d}'} = m + |\diag(i-1)|_{\mathfrak{d}} + |\diag(i+1)|_{\mathfrak{d}} - |\diag(i)|_{\mathfrak{d}} \]
\begin{proof}

Let $b$ denote the box of $\mathfrak{d}'$ with cylindric inversion coordinates $(i,i+1,0)$.
The sum of all the labels in boxes lying in the same cylindric column as $b$ is given by: 
\[ |\diag(i+1)|_{\mathfrak{d}} - |\diag(i)|_{\mathfrak{d}} \]
while the sum of all the labels in boxes lying strictly to the same cylindric row as $b$ is given by: 
\[ |\diag(i-1)|_{\mathfrak{d}} - |\diag(i)|_{\mathfrak{d}}\]
In otherwords:
\begin{align*}
\cohook(i,i+1,0)_{\mathfrak{d}'} & = m + |\diag(i-1)|_{\mathfrak{d}}  + |\diag(i+1)|_{\mathfrak{d}} - 2 |\diag(i)|_{\mathfrak{d}}
\end{align*}

Now:
\begin{align*}
|\diag(i)|_{\mathfrak{d}'} & = |\diag(i)|_{\mathfrak{d}} + \cohook(i,i+1,0)_{\mathfrak{d}'} \\
& = m + |\diag(i-1)|_{\mathfrak{d}} + |\diag(i+1)|_{\mathfrak{d}} - |\diag(i)|_{\mathfrak{d}}
\end{align*}
\end{proof}
\end{proposition}

\begin{proposition}\label{weight-preserving}
If $\psi_\pi$ is strongly weight preserving, then so is $\mathfrak{L}_i(\psi_\pi)$.
\begin{proof}

Let: 
\[ (m,\mathfrak{d}) = \mathfrak{l}_i[\mathfrak{d}'] \]
and let: 
\[\psi_\pi(\gamma, \mathfrak{d}) = \mathfrak{c} = (\mu^0, \ldots, \alpha, \nu, \beta, \ldots \mu^T)\] 

Then: 
\[\mathfrak{L}_i[\psi_\pi](\gamma, \mathfrak{d}') = \mathfrak{c}' = (\mu^0, \ldots, \alpha, \lambda, \beta, \ldots \mu^T)\]
where:
\[ \lambda = \mathfrak{U}_{\alpha,\beta}(m,\nu) \]

The weight condition for the local rule assures us that:
\[ |\lambda / \nu| = |\alpha / \nu| + |\beta / \nu| + m \]
In other words:
\[ |\lambda| = m + |\alpha| + |\beta| - |\nu| \]

For all $j \neq i$ we have 
\[|\diag(j)|_\mathfrak{d} = |\diag(j)|_{\mathfrak{d}'} = |\mu^j| - |\gamma| \]
We must show that:
\[ |\diag(i)|_{\mathfrak{d}'} = |\lambda| - |\gamma| \]

Now, by proposition \ref{weight-preserving} and the assumption that $\psi_\pi$ is strongly weight-preserving, we have:
\begin{align*}
|\diag(i)|_{\mathfrak{d}'} & = m +  |\diag(i-1)|_{\mathfrak{d}} + |\diag(i+1)|_{\mathfrak{d}} - |\diag(i)|_{\mathfrak{d}}\\
& = m + (|\alpha| - |\gamma|) + (|\beta| - |\gamma|) - (|\nu| - |\gamma|) \\
& = |\lambda| - |\gamma|
\end{align*}
The result follows.
\end{proof}
\end{proposition}

\newpage

\section{Conclusion}
We have shown that Fomin's growth diagram framework extends to the cylindric plane partition case. It remains to find an equivalent of the Hilmann-Grassl algorithm for the cylindric case. Other interesting questions include an investigation of the cyclic sieving phenomenon. Also, in the limiting case, is there an analog of the arctic circle theorem? Finally, Tingley \cite{tingley} showed that there is a bijection between cylindric plane partitions and crystal bases from the affine general linear group. It would be nice to understand what is the representation theoretic significance of Fomin's local rules.

\nocite{cylindric1}
\nocite{vanleeuwen}

\bibliographystyle{alpha}
\bibliography{references}

\end{document}

%% file: tableau.tex
\newdimen{\cellsize}
\newcommand\Bigboxes{\setlength{\cellsize}{24pt}\def\boxformat{}}%
\newcommand\bigboxes{\setlength{\cellsize}{18pt}\def\boxformat{}}
\newcommand\medboxes{\setlength{\cellsize}{14pt}\def\boxformat{}}
\newcommand\smallboxes{\setlength{\cellsize}{8pt}\def\boxformat{\scriptstyle}}
\medboxes
\newsavebox{\cellcontent}
\def\hidehrule#1#2{\kern-#1
  \hrule height#1 depth#2 \kern-#2 }%
\def\hidevrule#1#2{\kern-#1{\dimen\cellcontent=#1%
    \advance\dimen\cellcontent by#2\vrule width\dimen\cellcontent}\kern-#2 }%
\def\makeblankbox#1#2{\hbox{\lower\dp\cellcontent\vbox{\hidehrule{#1}{#2}%
    \kern-#1 
    \hbox to \wd\cellcontent{\hidevrule{#1}{#2}%
      \raise\ht\cellcontent\vbox to #1{}
      \lower\dp\cellcontent\vtop to #1{}
      \hfil\hidevrule{#2}{#1}}%
    \kern-#1\hidehrule{#2}{#1}}}}
\newcommand\cellify[1]{\defaultcell%
\sbox{\cellcontent}{\vbox to \cellsize{%
\vfill%
\hbox to \cellsize{\hfill$\boxformat #1$\hfill}
\vfill}}%
\rlap{\drawnbox}
\usebox{\cellcontent}}
\newcommand\tableau[1]{\vtop{\let\\\cr
\baselineskip -16000pt \lineskiplimit 16000pt \lineskip 0pt
\ialign{&\cellify{##}\cr#1\crcr}}}
\newcommand\defaultcell{\gdef\drawnbox{
\makeblankbox{0.2pt}{0.2pt}
}}
\newcommand\graycell{\gdef\drawnbox{%
\rlap{\color{Gray}\vrule width \cellsize height \cellsize}%
\makeblankbox{0.2pt}{0.2pt}
}}
\newcommand\bluecell{\gdef\drawnbox{%
\rlap{\color{Blue}\vrule width \cellsize height \cellsize}%
\makeblankbox{0.2pt}{0.2pt}
}}
\newcommand\redcell{\gdef\drawnbox{%
\rlap{\color{Red}\vrule width \cellsize height \cellsize}%
\makeblankbox{0.2pt}{0.2pt}
}}
\newcommand\greencell{\gdef\drawnbox{%
\rlap{\color{Green}\vrule width \cellsize height \cellsize}%
\makeblankbox{0.2pt}{0.2pt}
}}
\newcommand\orangecell{\gdef\drawnbox{%
\rlap{\color{Orange}\vrule width \cellsize height \cellsize}%
\makeblankbox{0.2pt}{0.2pt}
}}
\newcommand\yellowcell{\gdef\drawnbox{%
\rlap{\color{Yellow}\vrule width \cellsize height \cellsize}%
\makeblankbox{0.2pt}{0.2pt}
}}
\newcommand\thickyellowcell{\gdef\drawnbox{%
\rlap{\color{Yellow}\vrule width \cellsize height \cellsize}%
\makeblankbox{0.2pt}{0.1\cellsize}
}}
\newcommand\thickbluecell{\gdef\drawnbox{%
\rlap{\color{Blue}\vrule width \cellsize height \cellsize}%
\makeblankbox{0.2pt}{0.1\cellsize}
}}
\newcommand\thickredcell{\gdef\drawnbox{%
\rlap{\color{Red}\vrule width \cellsize height \cellsize}%
\makeblankbox{0.2pt}{0.1\cellsize}
}}
\newcommand\thickgreencell{\gdef\drawnbox{%
\rlap{\color{Green}\vrule width \cellsize height \cellsize}%
\makeblankbox{0.2pt}{0.1\cellsize}
}}
\newcommand\thickcell{\gdef\drawnbox{
\makeblankbox{0.2pt}{0.1\cellsize}%
}}
\newcommand\missingcell{\gdef\drawnbox{}}
\newcommand\vdotscell{\gdef\drawnbox{\kern-1.6pt\vbox{\baselineskip=4pt\lineskiplimit=0pt\hbox{}\hbox{.}\hbox{.}\hbox{.}\hbox{}}}}
\newcommand\hdotscell{\gdef\drawnbox{\vbox to \cellsize{\hbox{\kern1pt$\ldotp\ldotp\ldotp$}}}}
\newcommand\vhdotscell{\gdef\drawnbox{\rlap{\kern-1.6pt\vbox{\baselineskip=4pt\lineskiplimit=0pt\hbox{}\hbox{.}\hbox{.}\hbox{.}\hbox{}}}\vbox to \cellsize{\hbox{\kern1pt$\ldotp\ldotp\ldotp$}}}}
%
\newcommand\vertlinecell{\gdef\drawnbox{\unitlength=\cellsize%
\begin{picture}(1,1)
\put(0,0){\line(0,1){1}}
\end{picture}}}
\newcommand\horizlinecell{\gdef\drawnbox{\unitlength=\cellsize%
\begin{picture}(1,1)
\put(0,1){\line(1,0){1}}
\end{picture}}}
\newcommand\defaultcella{\gdef\drawnbox{
\unitlength=\cellsize%
\begin{picture}(1,1)
\put(0,0){\line(1,0){1}}
\put(0,0){\line(0,1){1}}
\put(1,0){\line(0,1){1}}
\put(0,1){\line(1,0){1}}
\end{picture}}}
\newcommand\thickcella{
\gdef\drawnbox{%
\unitlength=\cellsize%
\begin{picture}(1,1)
\linethickness{0.1\cellsize}
\put(0.0,0.05){\line(1,0){1}}
\put(0.05,0){\line(0,1){1}}
\put(0.95,0){\line(0,1){1}}
\put(0,0.96){\line(1,0){1}}
\end{picture}}}
\newcommand\graycella{\gdef\drawnbox{
\rlap{\color{Gray}\vrule width \cellsize height \cellsize}%
\unitlength=\cellsize%
\begin{picture}(1,1)
\put(0,0){\line(1,0){1}}
\put(0,0){\line(0,1){1}}
\put(1,0){\line(0,1){1}}
\put(0,1){\line(1,0){1}}
\end{picture}}%
}

%% file: article2.bbl
\begin{thebibliography}{GHT99}

\bibitem[Ada08]{adachi}
Shingo Adachi.
\newblock Reverse plane partitions and growth diagrams.
\newblock {\em preprint}, 2008.

\bibitem[Bor07]{borodin}
Alexei Borodin.
\newblock Periodic {S}chur process and cylindric partitions.
\newblock {\em Duke Math. J.}, 140(3):391--468, 2007.

\bibitem[Bur74]{burge}
W.H Burge.
\newblock Four correspondences between graphs and generalized {Y}oung tableau.
\newblock {\em J. Combin. Theory, Ser. A}, 17:12--30, 1974.

\bibitem[Fom86]{fomin1}
S.~V. Fomin.
\newblock The generalized {R}obinson-{S}chensted-{K}nuth correspondence.
\newblock {\em Zap. Nauchn. Sem. Leningrad. Otdel. Mat. Inst. Steklov. (LOMI)},
  155(Differentsialnaya Geometriya, Gruppy Li i Mekh. VIII):156--175, 195,
  1986.

\bibitem[Fom95]{fomin2}
Sergey Fomin.
\newblock Schur operators and {K}nuth correspondences.
\newblock {\em J. Combin. Theory Ser. A}, 72(2):277--292, 1995.

\bibitem[GHT99]{garsia}
A.~M. Garsia, M.~Haiman, and G.~Tesler.
\newblock Explicit plethystic formulas for {M}acdonald {$q,t$}-{K}ostka
  coefficients.
\newblock {\em S\'em. Lothar. Combin.}, 42:Art. B42m, 45 pp. (electronic),
  1999.
\newblock The Andrews Festschrift (Maratea, 1998).

\bibitem[GK97]{gessel-1997}
Ira~M. Gessel and C.~Krattenthaler.
\newblock Cylindric partitions.
\newblock {\em Trans. Amer. Math. Soc.}, 349(2):429--479, 1997.

\bibitem[Kra06]{krattenthaler}
C.~Krattenthaler.
\newblock Growth diagrams, and increasing and decreasing chains in fillings of
  ferrers shapes.
\newblock {\em Advances in Applied Mathematics}, 37:404--431, 2006.

\bibitem[Lan12]{cylindric1}
Robin Langer.
\newblock Enumeration of cylindric plane partitions {I}.
\newblock {\em available on arxiv}, 2012.

\bibitem[Mac95]{macdonald}
I.~G. Macdonald.
\newblock {\em Symmetric functions and {H}all polynomials}.
\newblock Oxford Mathematical Monographs. The Clarendon Press Oxford University
  Press, New York, second edition, 1995.
\newblock With contributions by A. Zelevinsky, Oxford Science Publications.

\bibitem[Oka96]{haskell}
C~Okasaki.
\newblock {\em Fully functional data structures}.
\newblock Cambridge University Press, 1996.

\bibitem[OR03]{okounkov}
Andrei Okounkov and Nikolai Reshetikhin.
\newblock Correlation function of {S}chur process with application to local
  geometry of a random 3-dimensional {Y}oung diagram.
\newblock {\em J. Amer. Math. Soc.}, 16(3):581--603 (electronic), 2003.

\bibitem[Tin08]{tingley}
Peter Tingley.
\newblock Three combinatorial models for {$\widehat{\rm sl}_n$} crystals, with
  applications to cylindric plane partitions.
\newblock {\em Int. Math. Res. Not. IMRN}, (2):Art. ID rnm143, 40, 2008.

\bibitem[vL05]{vanleeuwen}
Marc A.~A. van Leeuwen.
\newblock Spin-preserving {K}nuth correspondences for ribbon tableaux.
\newblock {\em Electron. J. Combin.}, 12:Research Paper 10, 65 pp.
  (electronic), 2005.

\end{thebibliography}
